\documentclass[preprint,review,12pt]{elsarticle}

\biboptions{sort&compress}
\usepackage[margin=3.5cm]{geometry}

\newcommand{\Lim}[1]{\raisebox{0.5ex}{\scalebox{0.8}{$\displaystyle \lim_{#1}\;$}}}

\usepackage{amssymb}
\usepackage{amsmath} 
\usepackage{amsfonts}
\usepackage{amsthm}
\usepackage{mathrsfs}
\newcommand{\norm}[1]{\left\lVert#1\right\rVert}
\usepackage{mathtools}
\usepackage{flushend}
\usepackage{upgreek}
\usepackage{yfonts}
\usepackage{algorithm}
\usepackage{algorithmic}

\journal{Elsevier Signal Processing}

\begin{document}

\begin{frontmatter}



\title{Private Networked Federated Learning\\ for Nonsmooth Objectives}


\author[inst1]{François Gauthier}
\author[inst1]{Cristiano Gratton}
\author[inst2]{Naveen K. D. Venkategowda}
\author[inst1]{Stefan Werner}

\affiliation[inst1]{organization={Department of Electronic Systems, NTNU},
            country={Norway}}

\affiliation[inst2]{organization={Department of Science and Technology, Linköping University},
            country={Sweden}}
\vspace{-10pt}
\begin{abstract}
This paper develops a networked federated learning algorithm to solve nonsmooth objective functions. 
To guarantee the confidentiality of the participants with respect to each other and potential eavesdroppers, we use the zero-concentrated differential privacy notion (zCDP). Privacy is achieved by perturbing the outcome of the computation at each client with a variance-decreasing Gaussian noise. ZCDP allows for better accuracy than the conventional $(\epsilon, \delta)$-DP and stronger guarantees than the more recent R\'enyi-DP by assuming adversaries aggregate all the exchanged messages.
The proposed algorithm relies on the distributed Alternating Direction Method of Multipliers (ADMM) and uses the approximation of the augmented Lagrangian to handle nonsmooth objective functions.
The developed private networked federated learning algorithm has a competitive privacy accuracy trade-off and handles nonsmooth and non-strongly convex problems.
We provide complete theoretical proof for the privacy guarantees and the algorithm's convergence to the exact solution. We also prove under additional assumptions that the algorithm converges in $O(1/n)$ ADMM iterations.
Finally, we observe the performance of the algorithm in a series of numerical simulations.

\end{abstract}



\vspace{-10pt}
\begin{keyword}
Federated Learning \sep Networked Architecture \sep Differential Privacy \sep Nonsmooth Objective Functions.
\end{keyword}

\end{frontmatter}


\section{Introduction}

Federated learning (FL) \cite{FLorigin} has garnered significant research attention recently because of its capacity to process massive amounts of data over a network of clients \cite{FedOverview, FLsurvey}. It has many applications, such as smart healthcare \cite{Healthcare}, autonomous vehicles \cite{AutonomousVehicles} and drones \cite{bedi2019asynchronous}, and industrial engineering \cite{FLapplications}. Networked FL \cite{NFL}, also called decentralized FL, proposes a collaborative approach to FL in which the problem is decomposed into many sub-problems that network clients solve by interacting with their immediate neighbors in a peer-to-peer fashion without involving a central coordinator. \cite{NFL_1, NFL_2, NFLex}. Networked FL is receiving growing interest as it resolves the limitations of using a single server in FL, such as communication and computation bottlenecks, while maintaining its advantages \cite{NFLsurvey}.

In many applications, the data held by clients is sensitive, and adversaries may try to extract private information from the information exchanged between the clients in the network. Therefore, it is imperative to mitigate information leakage during the client-interaction process in FL \cite{FLPrivacy}. In this context, differential privacy (DP) \cite{Dwork2006} provides a mechanism that protects individual privacy by ensuring minimal changes in the algorithm output, regardless of whether an individual data sample is present during the computation \cite{Dwork2014, OverviewDP}. Local DP protects the data of the clients with respect to each other and external eavesdroppers. This presents the advantage of protecting clients from honest-but-curious clients who form part of the network but share the information available to them with a third party. However, achieving good accuracy while providing high privacy guarantees in privacy-preserving FL is challenging, especially when several messages are exchanged.

To meet the demand for better privacy accuracy trade-off, concentrated differential privacy (CDP) was introduced in \cite{DworkCDP} as a relaxation of the conventional $(\epsilon,\delta)$-DP. Its purpose is to enable higher accuracy while maintaining identical privacy protection under the assumption that an adversary may aggregate all the exchanged messages \cite{jayaraman2019evaluating, Dwork2019}. CDP was relaxed into zero-concentrated differential privacy (zCDP) in \cite{Bun}, which is easier to use and offers similar benefits. Finally, dynamic-DP, introduced in \cite{zhang2016dual} can be used with zCDP to better suit iterative processes such as an ADMM algorithm. It enables iteration-specific privacy budgets. 
More recently, CDP has been relaxed into R\'enyi-DP \cite{RenyiDP}, which concentrates on a single moment of a privacy loss variable.
In contrast, CDP and zCDP provide a linear bound on all positive moments and, therefore, a stronger privacy guarantee. In this work, we use dynamic zCDP.

Existing networked FL and distributed learning solutions mainly comprise (sub)gradient-based and ADMM-based algorithms. The former typically converge at a rate of $O(1/\sqrt{n})$ \cite{Nedic2009}, and the latter usually converge at a rate of $O(1/n)$ \cite{Zhuhan}. Although networked FL is recent, privacy-preserving distributed learning has been thoroughly studied \cite{Zhuhan,Zhang2017,Zhang2018,ESPPrivADMM,ding2019differentially,zhang2016dual,Renyi}.
Among these works, we can classify two relevant groups to the problem at hand. The first group comprises solutions that use a networked architecture and assume the objective functions to be smooth and convex \cite{Zhuhan,Zhang2017,Zhang2018,ESPPrivADMM,ding2019differentially,zhang2016dual},
which may not be a valid assumption in practice. In fact, many compelling objectives cannot be accurately modeled in this way \cite{gentle1977least,roth2004generalized}. Among the above works, \cite{Zhuhan} offers a convergence rate of $O(1/n)$ and, in \cite{Renyi}, the regularizer function can be nonsmooth, but the loss function is assumed smooth and differentiable. In addition, the convergence rate of the algorithm in \cite{Renyi} is $O(1/\sqrt{n})$. The second group comprises solutions relying on a single server and handling nonsmooth and non-strongly convex objectives \cite{NSLc1,NSLc2,NSLc3,Huang2020,hu2019learning}. We note that the works in \cite{NSLc1,NSLc2,NSLc3} do not consider privacy, and the work in \cite{hu2019learning} can only accommodate nonsmooth regularizer functions. The solution proposed in \cite{Huang2020} handles nonsmooth and non-strongly convex objectives but converges in $O(1/\sqrt{n})$ and relies on a server to aggregate the local models.

This paper proposes a privacy-preserving networked FL algorithm that handles nonsmooth and non-strongly convex objective functions. Furthermore, the proposed algorithm is proven to converge to the exact solution with a rate of $O(1/n)$.
We consider a distributed network of clients that solve an optimization problem collaboratively. Each client iteratively updates its local model using its local data and the models received from its neighbors. To ensure the confidentiality of its local data, a client perturbs its model before communication with white Gaussian noise. The variance of the noise added to the models is such that the total privacy leakage of the clients throughout the computation is bounded under the zCDP metric. The ADMM used to solve the optimization problem is distributed to comply with the networked setting. Further, the clients' primal updates use an approximation of the augmented Lagrangian obtained by taking the first-order approximation of the objective function. This enables the proposed method to handle nonsmooth objective functions. Mathematical analysis shows that the proposed algorithm converges to the optimal point in $O(1/n)$. Numerical simulations compare the proposed algorithm with the most closely related distributed learning solutions.

The rest of the manuscript is organized as follows. Section II introduces the proposed zero-Concentrated Differentially Private Networked Federated Learning (zCDP-NFL) algorithm. Sections III and IV contain the mathematical analysis of the privacy guarantees and convergence properties of the proposed algorithm, respectively. Section V compares zCDP-NFL with existing methods in a series of numerical simulations.

\textit{Mathematical notations:} Matrices, column vectors, and scalars will be respectively denoted by bold uppercase, bold lowercase, and lowercase letters. 
The set of natural integers is denoted by $\mathbb{N}$ and the set of real numbers by $\mathbb{R}$. 
The operators $(\cdot)^\mathsf{T}$ denotes the transpose of a matrix. 
$|| \cdot ||$ represents the Euclidean norm, $|| \cdot ||_1$ the $L_1$ norm, and $||\boldsymbol{x}||^2_{\boldsymbol{G}} = \langle\boldsymbol{x},\boldsymbol{G}\boldsymbol{x}\rangle$ for any couple of vector $\boldsymbol{x}$ and matrix $\boldsymbol{G}$. 
The inner product between two vectors $\boldsymbol{a}$ and $\boldsymbol{b}$ is denoted by  $\langle \boldsymbol{a},\boldsymbol{b} \rangle$. 
The statistical expectation operator is represented by $\mathbb{E}[\cdot]$ and $\mathcal{N}(\boldsymbol{\mu},\boldsymbol{\Sigma})$ denotes the normal distribution with mean $\boldsymbol{\mu}$ and covariance matrix $\boldsymbol{\Sigma}$. 
If a random variable $A$ follows the law $\mathcal{B}$, we will write $A \sim \mathcal{B}$. The identity matrix in $\mathbb{R}^{n \times n}$ is denoted by $\boldsymbol{I}_n$ and subgradient of a function $g(\cdot)$ is denoted by $g'(\cdot)$. 
The nonzero smallest and largest singular values of a semidefinite matrix $\boldsymbol{A}$ are denoted by  $\Phi_{\min}(\boldsymbol{A})$ and $\Phi_{\max}(\boldsymbol{A})$.
\vspace{-10pt}
\section{Approximated Private Networked Federated Learning}


\subsection{Distributed Empirical Risk Minimization}

We consider a connected network of clients modeled as an undirected graph $\mathcal{G}(\mathcal{C},\mathcal{E})$ where vertex set $\mathcal{C}= \{1, \hdots, K\}$ corresponds to the clients and edge set $\mathcal{E}$ contains the $|\mathcal{E}| = E$ undirected communication links. The set $\mathcal{N}_{k}$ contains the indexes of the neighbors of client $k$.

Each client $k\in\mathcal{C}$ has a private data set $\mathcal{D}_k \coloneqq\{(\mathbf{X}_k,\mathbf{y}_k):\mathbf{X}_k=[\mathbf{x}_{k,1},\hdots,\mathbf{x}_{k,M_k}]^\mathsf{T} \in\mathbb{R}^{M_k\times P}, \ \mathbf{y}_k=[y_{k,1},\hdots,y_{k,M_k}]^\mathsf{T}\in\mathbb{R}^{M_k}\}$, where $M_k$ is the number of data samples and $P$ the number of features in the data.


To fit with the networked architecture, we consider the distributed empirical risk minimization problem with local primal variables $\mathcal{V}\coloneqq\{{\bf w}_k\}_{k=1}^K$:
\begin{equation}
\begin{aligned}
&\underset{\{{\bf w}_k\}}{\min}
&& \sum_{k=1}^{K}\Bigl(\frac{1}{M_k} \sum_{j=1}^{M_k} \ell(\mathbf{x}_{k,j},\mathbf{y}_{k,j};{\bf w}_k)+\frac{\lambda}{K} R({\bf w}_k) \Bigr) \\
&\text{\ s.t.}\ 
&& {\bf w}_k=\mathbf{z}_k^l,\ {\bf w}_l=\mathbf{z}_k^l, \quad l\in\mathcal{N}_k, \quad\forall k\in\mathcal{C},
\end{aligned}
\label{distribERM}
\end{equation}
where $\ell:\mathbb{R}^{P}\rightarrow\mathbb{R}$ is the loss function, $R:\mathbb{R}^{P}\rightarrow\mathbb{R}$ is the regularizer function, $\lambda>0$ is the regularization parameter, and the equality constraints enforce consensus. The auxiliary variables $\mathcal{Z}\coloneqq\{\mathbf{z}_k^l\}_{l \in \mathcal{N}_k}$ are only used to derive the local recursions and are eventually eliminated. In the following, we consider the learning problem where $\ell(\cdot)$ and $R(\cdot)$ are convex, but not necessarily strongly convex or smooth.


\subsection{Approximate Augmented Lagrangian}

The augmented Lagrangian associated with \eqref{distribERM} is given by
\begin{align}
\label{Lagrangian}
    \mathcal{L}_{\rho}&(\mathcal{V},\mathcal{M},\mathcal{Z}) = \sum_{k=1}^K \Big( \frac{\ell(\mathbf{X}_k,\mathbf{y}_k;{\bf w}_k)}{M_k} + \frac{\lambda R({\bf w}_k)}{K} \Big) \notag \\
    &+ \sum_{k=1}^K \sum_{l \in \mathcal{N}_k} \Big[ \boldsymbol{\mu}_k^{l \mathsf{T}}({\bf w}_k-\mathbf{z}_k^l) + \boldsymbol{\gamma}_k^{l \mathsf{T}}({\bf w}_l-\mathbf{z}_k^l) \Big] \\
    &+ \frac{\rho}{2} \sum_{k=1}^K \sum_{l \in \mathcal{N}_k} \Big( ||{\bf w}_k-\mathbf{z}_k^l||^2 + ||{\bf w}_l-\mathbf{z}_k^l||^2 \Big) \notag
\end{align}
where $\rho>0$ is a penalty parameter and $\mathcal{M}\coloneqq\{\{\boldsymbol{\mu}_k^l\}_{l\in\mathcal{N}_k},\{\boldsymbol{\gamma}_k^l\}_{l\in\mathcal{N}_k}\}_{k=1}^K$ are the Lagrange multipliers associated with the constraints in \eqref{distribERM}.

Given that the Lagrange multipliers $\mathcal{M}$ are initialized to zero, by using the Karush-Kuhn-Tucker conditions of optimality for \eqref{distribERM} and setting $\boldsymbol{\gamma}_k^{(n)}=2\sum_{l\in\mathcal{N}_k}(\boldsymbol{\gamma}_k^l)^{(n)}$, it can be shown that the Lagrange multipliers $\{\boldsymbol{\mu}_k^l\}_{l\in\mathcal{N}_k}$ and the auxiliary variables $\mathcal{Z}$ are eliminated \cite{Giannakis2016,Dridge}. The resulting algorithm reduces to the following iterative steps at client $k$.
\begin{align}
{\bf w}_k^{(n)} &= \arg\min_{{\bf w}_k} \Bigg[ f_k({\bf w}_k) + {\bf w}_k^\mathsf{T}\boldsymbol{\gamma}_k^{(n-1)} + \rho\sum_{l\in\mathcal{N}_k}\bigg\|{\bf w}_k-\frac{{\bf w}_k^{(n-1)}+{\bf w}_l^{(n-1)}}{2}\bigg\|^2 \Bigg] \notag \\
    \boldsymbol{\gamma}_k^{(n)} &= \boldsymbol{\gamma}_k^{(n-1)} {+} \rho\sum_{l\in\mathcal{N}_k}\left({\bf w}_k^{(n)} - {\bf w}_l^{(n)}\right)
\end{align}
where $n$ is the iteration index and
\begin{align}
    f_k({\bf w}_k) = \frac{\ell(\mathbf{X}_k,\mathbf{y}_k;{\bf w}_k)}{M_k} + \frac{\lambda R({\bf w}_k)}{K}.
\end{align}

To handle nonsmooth $\ell(\cdot)$ and $R(\cdot)$ functions, we take the first-order approximation of $f_k$ with an $l_2$-norm prox function, denoted as $\hat{f_k}$. Similarly as in \cite{Huang2020,nemirovski2009robust}, such an approximation is given by
\begin{align}
\hat{f_k}({\bf w}_k;\mathcal{V}^{(n)}) &= \frac{\ell(\mathbf{X}_k,\mathbf{y}_k;{\bf w}_k^{(n)})}{M_k} + \frac{\lambda R({\bf w}_k^{(n)})}{K} + \frac{\big\|{\bf w}_k-{\bf w}_k^{(n)}\big\|^2}{2\eta_k^{(n+1)}} \label{objective} \\
&+\left({\bf w}_k-{\bf w}_k^{(n)}\right)^{\mathsf{T}}\Bigl(\frac{\ell'(\mathbf{X}_k,\mathbf{y}_k;{\bf w}_k^{(n)})}{M_k} + \frac{\lambda R'({\bf w}_k^{(n)})}{K}\Bigr) \notag
\end{align}
where $\mathcal{V}^{(n)}=\{ {\bf w}_k^{(n)}, k \in \mathcal{C} \}$, and $\eta_k^{(n)}$ is a time-varying step size.

Taking the first-order approximation of $f_k$ leads to an inexact update at a given iteration; however, the algorithm does not need to solve the problem with high precision at each iteration to guarantee overall accuracy \cite{Huang2020}. In the end, considering $\hat{f}_k$ instead of $f_k$ in the primal update makes the algorithm capable of solving nonsmooth objectives with a minimal impact on overall accuracy. Unlike the method used in \cite{Nedic2009} to deal with nonsmooth objective functions, the approach taken here is compatible with the algorithm's convergence to the exact objective value.

\subsection{Privacy Preservation}

To prevent the leakage of private information, we introduce local differential privacy to the algorithm via message perturbation. For this purpose, each client $k$ shares at iteration $n$ with its neighbors the perturbed estimate
\begin{align}
    \widetilde{{\bf w}}_k^{(n)} & =   {\bf w}_k^{(n)} {+} \boldsymbol{\xi}_k^{(n)}  \label{update2}
\end{align}
with $\boldsymbol{\xi}_k^{(n)} \sim \mathcal{N}(\mathbf{0},\sigma_{k}^{2}(n)\mathbf{I}_P)$. We denote $\widetilde{\mathcal{V}}^{(n)}=\{\widetilde{{\bf w}}_k^{(n)}, k \in \mathcal{C} \}$.

The value of the noise perturbation variance, $\sigma_{k}^{2}(n)$, in \eqref{update2} dictates the privacy protection of the algorithm.
To guarantee convergence to the optimal solution, as opposed to a neighborhood of it, the variance must decrease with the iterations \cite{huang2015differentially}. This is made possible by using dynamic zCDP, where the privacy budget is iteration-specific.

The proposed zero-Concentrated Differentially Private Networked Federated Learning (zCDP-NFL) algorithm is developed in algorithm \ref{alg_zCDP_NFL}. It is a networked federated learning algorithm that protects the privacy of the clients with local dynamic differential privacy and can handle nonsmooth objective functions. In the following sections, we conduct mathematical analysis to quantify its privacy protection and establish convergence guarantees.
\begin{algorithm}[h!]
\caption{zCDP-NFL}
\label{alg_zCDP_NFL}
\begin{algorithmic}[1]
    \item[] \textbf{Initialization:} ${\bf w}_k^{(0)}=\mathbf{0}$, $\boldsymbol{\gamma}_k^{(0)}=\mathbf{0}$, $\forall k \in \mathcal{K}$
    \item[] \textit{-- Procedure at client $k$ --}
    \item[] \textbf{For} iteration $n = 1, 2, \hdots$:
    \begin{align}
        \label{PrimalUpdate}
        {\bf w}_k^{(n)} &= \arg\min_{{\bf w}_k}   \hat{f_k}({\bf w}_k;\widetilde{\mathcal{V}}^{(n-1)}) + {\bf w}_k^\mathsf{T}\boldsymbol{\gamma}_k^{(n-1)} + \rho\sum_{l\in\mathcal{N}_k}\bigg\|{\bf w}_k-\frac{\widetilde{{\bf w}}_k^{(n-1)}+\widetilde{{\bf w}}_l^{(n-1)}}{2}\bigg\|^2  \\
        \label{NoisePerturbation}
        \widetilde{{\bf w}}_k^{(n)} &=   {\bf w}_k^{(n)} {+} \boldsymbol{\xi}_k^{(n)} \\
        \label{DualUpdate}
        \boldsymbol{\gamma}_k^{(n)} & = \boldsymbol{\gamma}_k^{(n-1)} {+} \rho\sum_{l\in\mathcal{N}_k}\left(\widetilde{{\bf w}}_k^{(n)} {-} \widetilde{{\bf w}}_l^{(n)}\right)
    \end{align}
    \item[] \textbf{End For}
\end{algorithmic} 
\end{algorithm}

\section{Privacy Analysis}

The first step in the quantification of differential privacy is to measure the impact of an individual data sample on the output of the local training process. For this purpose, we define the $l_2$-norm sensitivity as follows.


\noindent \textbf{Definition I.} \textit{The $l_2$-norm sensitivity is given by}
\begin{equation}
\Delta_{k,2}=\max_{\mathcal{D}_k,\mathcal{D}_k'}\norm{{\bf w}_{k,\mathcal{D}_k}^{(n)}-{\bf w}_{k,\mathcal{D}_k'}^{(n)}}
\end{equation}
\textit{where ${\bf w}_{k,\mathcal{D}_k}^{(n)}$ and ${\bf w}_{k,\mathcal{D}_k'}^{(n)}$ denote the local primal variable updates from two neighboring data sets $\mathcal{D}_k$ and $\mathcal{D}_k'$ differing in only one data sample $(\mathbf{x}_{k,M_k}',y_{k,M_k}')$, i.e., $\mathcal{D}_k'\coloneqq\{(\mathbf{X}_k',\mathbf{y}_k'):\mathbf{X}_k'=[\mathbf{x}_{k,1},\mathbf{x}_{k,2},\hdots,\mathbf{x}_{k,M_k-1},\mathbf{x}_{k,M_k}']^\mathsf{T}\in\mathbb{R}^{M_k\times P}, \ \mathbf{x}_{k,j}\in\mathbb{R}^P, \ j=1,\hdots,M_k, \mathbf{y'}_k=[y_{k,1},y_{k,2},\hdots,y_{k,M_k-1},y_{k,M_k}']^\mathsf{T}\in\mathbb{R}^{M_k}\}$.}

Two parameters govern privacy protection in dynamic zCDP. The initial privacy value, $\varphi_{k}^{(0)}$, and the variance decrease rate, $\tau$. To establish a relation between the privacy value at a given iteration, $\varphi_{k}^{(n)}$, and the noise perturbation, it is necessary to take the following assumption.


\noindent \textbf{Assumption 1.} \textit{The functions $\ell_k(\cdot)$ have bounded gradient, that is, there exists a constant $c_1$ such that $||\ell_k'(\cdot)|| \leqslant c_1, \forall k \in \mathcal{C}$.}


We now quantify the $l_2$-norm sensitivity in the following result.

\noindent \textbf{Lemma I.} \textit{Under Assumption 1, the $l_2$-norm sensitivity is given by
\begin{equation}
    \Delta_{k,2}(n) = \max_{\mathcal{D},\mathcal{D}'} ||{\bf w}_{k,\mathcal{D}}^{(n)} - {\bf w}_{k,\mathcal{D}'}^{(n)}|| \leqslant \frac{2 c_1}{M_k (2 \rho |\mathcal{N}_k| + \frac{1}{\eta^{(n)}})}.
\end{equation}}
\begin{proof}
See Appendix A.
\end{proof}

With the $l_2$-norm sensitivity, we can establish the relation between the noise perturbation in \eqref{NoisePerturbation} and the privacy value $\varphi_{k}^{(n)}$, quantifying the local privacy guarantee of the algorithm in terms of zCDP.

\noindent \textbf{Theorem I.} \textit{Under Assumption 1, zCDP-NFL satisfies dynamic $\varphi_{k}^{(n)}$-zCDP with the relation between $\varphi_{k}^{(n)}$ and $\sigma_{k}^{2}(n)$ given by}
\begin{equation}
\label{secondeqdeftwo}
\sigma_{k}^{2}(n)=\frac{\Delta_{k,2}^{2}(n)}{2\varphi_{k}^{(n)}}.
\end{equation}
\begin{proof}
See Appendix B.
\end{proof}

Using the result above, it is possible to obtain the total privacy guarantee throughout the computation in terms of $(\epsilon,\delta)$-DP using \cite[Lemma 1.7]{Bun}. We establish the following.

\noindent \textbf{Corollary.} \textit{For any $\tau \in (0,1)$ and $\delta \in (0,1)$, zCDP-NFL guarantees $(\epsilon,\delta)$-DP throughout the computation with $\epsilon = \displaystyle\max_{k \in \mathcal{C}} \epsilon_k$, where $\epsilon_k = \varphi_{k}^{(1)} \frac{1 - \tau^{T}}{\tau^{T-1} - \tau^{T}} + 2 \sqrt{\varphi_{k}^{(1)} \frac{1 - \tau^{T}}{\tau^{T-1} - \tau^{T}} \log \frac{1}{\delta}}$, and $T$ is the final iteration index.}
\begin{proof}
See Appendix C.
\end{proof}


\section{Convergence Analysis}
This section proves that the zCDP-NFL algorithm converges to the optimal value in $O(1/n)$ iterations under the following assumption that the objective function $f(\cdot)$ is convex. Additionally, we derive the privacy-accuracy trade-off bound of the algorithm.

\noindent \textbf{Assumption 2.} \textit{The objective function $f(\cdot)$ is convex.}

\subsection{Alternative Representation}
We begin by transforming the minimization problem \eqref{distribERM} into \eqref{ADMM2} by reformulating the conditions. We denote by $\boldsymbol{\mathfrak{w}}=[{\bf w}_1^\mathsf{T},{\bf w}_2^\mathsf{T},...,{\bf w}_K^\mathsf{T}]^\mathsf{T} \in \mathbb{R}^{K P}$, and $\boldsymbol{\mathfrak{z}} = [(\boldsymbol{z}_k^l)^\mathsf{T}, (\boldsymbol{z}_l^k)^\mathsf{T}; \forall (k,l) \in \mathcal{E}]^\mathsf{T} \in \mathbb{R}^{2E P}$ the vectors of the concatenated vectors ${\bf w}_k$ and $\boldsymbol{z}_k^l$ respectively. We also introduce the matrices $\boldsymbol{A}_1, \boldsymbol{A}_2 \in \mathbb{R}^{2EP \times KP}$ composed of $P \times P$-sized blocks. Given a couple of connected clients $(k, l) \in \mathcal{E}$, their associated auxiliary variable $\boldsymbol{z}_{k,l}$, and its corresponding index in $\boldsymbol{\mathfrak{z}}$, $q$; the blocks $ \bigl( {\bf A}_1 \bigr)_{q, k} $ and $ \bigl( {\bf A}_2 \bigr)_{q, l} $ are equal to the identity matrix ${\bf I}_d$, all other blocks are null. Finally, we set $\boldsymbol{A}=[\boldsymbol{A}_1;\boldsymbol{A}_2] \in \mathbb{R}^{4 E P \times K P}$ and $\boldsymbol{B}=[-\boldsymbol{I}_{2EP};-\boldsymbol{I}_{2EP}]  \in \mathbb{R}^{4 E P \times 2 E P}$. Hence, we can reformulate \eqref{distribERM} as
%
\begin{equation}
\begin{aligned}
&\underset{\boldsymbol{\mathfrak{w}}}{\min}
&& \sum_{k=1}^{K}\Bigl(\frac{1}{M_k} \sum_{j=1}^{M_k} \ell(\mathbf{x}_{k,j},\mathbf{y}_{k,j};{\bf w}_k)+\frac{\lambda}{K} R({\bf w}_k) \Bigr). \\
&\text{\ s.t.}\ 
&&\boldsymbol{A} \boldsymbol{\mathfrak{w}} + \boldsymbol{B} \boldsymbol{\mathfrak{z}} = 0
\end{aligned}
\label{ADMM2}
\end{equation}

The newly introduced matrices can be used to reformulate the Lagrangian, the objective function, and the ADMM steps.
The conventional augmented Lagrangian in \eqref{Lagrangian} can be expressed as
\begin{align}
    \mathcal{L}_\rho = f(\boldsymbol{\mathfrak{w}},\widetilde{\mathcal{V}}^{(n)}) + \langle\boldsymbol{A} \boldsymbol{\mathfrak{w}} + \boldsymbol{B} \boldsymbol{\mathfrak{z}}, \boldsymbol{\lambda}\rangle + \frac{\rho}{2} ||\boldsymbol{A} \boldsymbol{\mathfrak{w}} + \boldsymbol{B} \boldsymbol{\mathfrak{z}}||^2 \notag
\end{align}
where $f(\boldsymbol{\mathfrak{w}},\widetilde{\mathcal{V}}^{(n)}) = \sum_{k=1}^K f({\bf w}_k,\widetilde{\mathcal{V}}^{(n)})$. Similarly, the augmented Lagrangian, corresponding to the use of the first-order approximation of the objective function in \eqref{objective}, can be expressed as
\begin{align}
    \hat{\mathcal{L}}_\rho = \hat{f}(\boldsymbol{\mathfrak{w}},\widetilde{\mathcal{V}}^{(n)}) + \langle\boldsymbol{A} \boldsymbol{\mathfrak{w}} + \boldsymbol{B} \boldsymbol{\mathfrak{z}}, \boldsymbol{\lambda}\rangle + \frac{\rho}{2} ||\boldsymbol{A} \boldsymbol{\mathfrak{w}} + \boldsymbol{B} \boldsymbol{\mathfrak{z}}||^2 \notag
\end{align}
where $\hat{f}(\boldsymbol{\mathfrak{w}},\widetilde{\mathcal{V}}^{(n)}) = \sum_{k=1}^K \hat{f}_k({\bf w}_k,\widetilde{\mathcal{V}}^{(n)})$ with $\hat{f}_k({\bf w}_k,\widetilde{\mathcal{V}}^{(n)})$, as defined in \eqref{objective}. 

From now on, we will denote $\hat{f}(\boldsymbol{\mathfrak{w}},\widetilde{\mathcal{V}}^{(n)})$ and $\hat{f}_k({\bf w}_k,\widetilde{\mathcal{V}}^{(n)})$ by $\hat{f}(\boldsymbol{\mathfrak{w}})$ and $\hat{f}_k({\bf w}_k)$, respectively. 
Further, we  let $\widetilde{\boldsymbol{\mathfrak{w}}}^{(n)}$, $\boldsymbol{\mathfrak{w}}^{(n)}$, and $\boldsymbol{\xi}^{(n)}$ denote the concatenation of $\widetilde{{\bf w}}_k^{(n)}$, ${\bf w}_k^{(n)}$, and $\boldsymbol{\xi}_k^{(n)}$, respectively, such that $\widetilde{\boldsymbol{\mathfrak{w}}}^{(n)} = \boldsymbol{\mathfrak{w}}^{(n)} + \boldsymbol{\xi}^{(n)}$.

We introduce the diagonal matrix $\boldsymbol{D}^{(n+1)} \in \mathbb{R}^{K \times K}$ comprising the time-varying step sizes, i.e., $[\boldsymbol{D}^{(n+1)}]_{k,k} = \frac{1}{\sqrt{2 \eta_k^{(n+1)}}}$, and reformulate $\hat{f}(\boldsymbol{\mathfrak{w}}^{(n+1)})$ in matrix form: \vspace{-10pt}
\begin{align}
\label{fhat_matrix}
    \hat{f}(\boldsymbol{\mathfrak{w}}^{(n+1)}) = & f(\widetilde{\boldsymbol{\mathfrak{w}}}^{(n)}) + ||\boldsymbol{D}^{(n+1)} \otimes \boldsymbol{I}_P (\boldsymbol{\mathfrak{w}}^{(n+1)} - \widetilde{\boldsymbol{\mathfrak{w}}}^{(n)}) ||^2 \notag\\
    &+ (\boldsymbol{\mathfrak{w}}^{(n+1)} - \widetilde{\boldsymbol{\mathfrak{w}}}^{(n)})^\mathsf{T} f'(\widetilde{\boldsymbol{\mathfrak{w}}}^{(n)})
\end{align}

\noindent The resulting function $\hat{f}$ is convex with respect to $\boldsymbol{\mathfrak{w}}$. That is, it satisfies $\hat{f}(\widetilde{\boldsymbol{\mathfrak{w}}}^{(n)}) - \hat{f}(\boldsymbol{\mathfrak{w}}) \leqslant \langle \widetilde{\boldsymbol{\mathfrak{w}}}^{(n)} - \boldsymbol{\mathfrak{w}}, \hat{f}'(\widetilde{\boldsymbol{\mathfrak{w}}}^{(n)})\rangle$, where the subgradient $\hat{f}'(\boldsymbol{\mathfrak{w}}^{(n+1)}) \in \mathbb{R}^{K P}$ is given by
    $\hat{f}'(\boldsymbol{\mathfrak{w}}^{(n+1)}) = 2 \boldsymbol{D}^{(n+1)} \otimes \boldsymbol{I}_P (\boldsymbol{\mathfrak{w}}^{(n+1)} - \widetilde{\boldsymbol{\mathfrak{w}}}^{(n)}) + f'(\widetilde{\boldsymbol{\mathfrak{w}}}^{(n)})$.

The steps of the ADMM, consisting of  the minimization of $\hat{\mathcal{L}}_\rho$ with respect to $\boldsymbol{\mathfrak{w}}, \boldsymbol{\mathfrak{z}}$ and $\boldsymbol{\lambda}$ alternatively, can now be reformulated with the newly introduced variables as follows:\vspace{-10pt}
\begin{align}
    \hat{f}'(\boldsymbol{\mathfrak{w}}^{(n+1)}) + \boldsymbol{A}^\mathsf{T} \boldsymbol{\lambda}^{(n)} + \rho \boldsymbol{A}^\mathsf{T} (\boldsymbol{A} \boldsymbol{\mathfrak{w}}^{(n+1)} + \boldsymbol{B} \boldsymbol{\mathfrak{z}}^{(n)}) = 0 \notag \\
    \boldsymbol{B}^\mathsf{T} \boldsymbol{\lambda}^{(n)} + \rho \boldsymbol{B}^\mathsf{T} (\boldsymbol{A} \widetilde{\boldsymbol{\mathfrak{w}}}^{(n+1)} + \boldsymbol{B} \boldsymbol{\mathfrak{z}}^{(n+1)}) = 0 \label{CVproofADMM1} \\
    \boldsymbol{\lambda}^{(n+1)} - \boldsymbol{\lambda}^{(n)} + \rho (\boldsymbol{A} \widetilde{\boldsymbol{\mathfrak{w}}}^{(n+1)} + \boldsymbol{B} \boldsymbol{\mathfrak{z}}^{(n+1)}) = 0 \notag
\end{align}
We introduce the following auxiliary matrices in order to reduce \eqref{CVproofADMM1} to two steps, similarly as in \cite{shi2014linear}: $\boldsymbol{H}_+ = \boldsymbol{A}_1^\mathsf{T} + \boldsymbol{A}_2^\mathsf{T}$, $\boldsymbol{H}_- = \boldsymbol{A}_1^\mathsf{T} - \boldsymbol{A}_2^\mathsf{T}$, $\boldsymbol{\alpha} = \boldsymbol{H}_-^\mathsf{T} \boldsymbol{\mathfrak{w}}$, $\boldsymbol{L}_+ = \frac{1}{2} \boldsymbol{H}_+ \boldsymbol{H}_+^\mathsf{T}$, $\boldsymbol{L}_- = \frac{1}{2} \boldsymbol{H}_- \boldsymbol{H}_-^\mathsf{T}$ and $\boldsymbol{M} = \frac{1}{2} (\boldsymbol{L}_+ + \boldsymbol{L}_-)$. We note that $\boldsymbol{L}_+$ and $\boldsymbol{L}_-$ correspond to the signless Laplacian and signed Laplacian matrices of the network, respectively. Hence, $\boldsymbol{L}_-$ is positive semi-definite with the nullspace given by $Null(\boldsymbol{L}_-) = span\{1\}$. Then, as derived in \cite[Section II.B]{shi2014linear}, \eqref{CVproofADMM1} becomes \vspace{-10pt}
\begin{align}\hat{f}'(\boldsymbol{\mathfrak{w}}^{(n+1)}) + \boldsymbol{\alpha}^{(n)} + 2 \rho \boldsymbol{M} \boldsymbol{\mathfrak{w}}^{(n+1)} - \rho \boldsymbol{L}_+ \widetilde{\boldsymbol{\mathfrak{w}}}^{(n)} = 0 \label{CVproofPrivateADMM} \\
    \boldsymbol{\alpha}^{(n+1)} - \boldsymbol{\alpha}^{(n)} - \rho \boldsymbol{L}_- \widetilde{\boldsymbol{\mathfrak{w}}}^{(n+1)} = 0 \notag
\end{align}
%
%
The last reformulation step is based on the work in \cite{li2017robust}. We introduce the matrix $\boldsymbol{Q} = \sqrt{\boldsymbol{L}_- / 2}$, note that by construction $Null(\boldsymbol{Q}) = span\{1\}$, the auxiliary sequence $\boldsymbol{r}^{(n)} = \sum_{s=0}^{n} \boldsymbol{Q} \widetilde{\boldsymbol{\mathfrak{w}}}^{(s)}$, vector $\boldsymbol{q}^{(n)} = \dbinom{\boldsymbol{r}^{(n)}}{\widetilde{\boldsymbol{\mathfrak{w}}}^{(n)}}$, and matrix $\boldsymbol{G} = \dbinom{\rho \boldsymbol{I} \; \; \; \; \; \; 0}{0 \; \; \rho \boldsymbol{L}_+ / 2}$. Combining both equations in \eqref{CVproofPrivateADMM}, as in \cite[Lemma 1]{li2017robust}, and reformulating the result, see see \cite[Lemma 2]{li2017robust}, we obtain
\begin{align}
    \frac{\hat{f}'(\boldsymbol{\mathfrak{w}}^{(n+1)})}{\rho} + 2 \boldsymbol{Q} \boldsymbol{r}^{(n+1)} + \boldsymbol{L}_+(\boldsymbol{\mathfrak{w}}^{(n+1)} - \widetilde{\boldsymbol{\mathfrak{w}}}^{(n)}) = 2 \boldsymbol{M} \boldsymbol{\xi}^{(n+1)}.
    \label{CVproofIntroLast}
\end{align}

\subsection{Convergence Proof}

We start by establishing a bound for the distance to the optimal solution, denoted $\boldsymbol{\mathfrak{w}}^{*}$, at a given iteration.

\noindent \textbf{Lemma II.} \textit{For any $\boldsymbol{r} \in \mathbb{R}^{K P}$ and at any iteration $n$, we have
\begin{align}
    &\frac{f(\widetilde{\boldsymbol{\mathfrak{w}}}^{(n)}) - f(\boldsymbol{\mathfrak{w}}^{*})}{\rho} + \langle\widetilde{\boldsymbol{\mathfrak{w}}}^{(n)}, 2 \boldsymbol{Q} \boldsymbol{r}\rangle \\
    &\leqslant \frac{1}{\rho}(||\boldsymbol{q}^{(n-1)} - \boldsymbol{q}^{*}||^2_{\boldsymbol{G}} - ||\boldsymbol{q}^{(n)} - \boldsymbol{q}^{*}||^2_{\boldsymbol{G}})   - 2 \langle \boldsymbol{Q} \widetilde{\boldsymbol{\mathfrak{w}}}^{(n)}, \boldsymbol{Q} \widetilde{\boldsymbol{\mathfrak{w}}}^{(n+1)} \rangle - ||\widetilde{\boldsymbol{\mathfrak{w}}}^{(n)} - \widetilde{\boldsymbol{\mathfrak{w}}}^{(n-1)}||_{\frac{\boldsymbol{L}_+}{2}}^2 \notag\\
    &\quad + \langle \widetilde{\boldsymbol{\mathfrak{w}}}^{(n)} - \boldsymbol{\mathfrak{w}}^{*}, \boldsymbol{L}_+(2 \widetilde{\boldsymbol{\mathfrak{w}}}^{(n)} - \widetilde{\boldsymbol{\mathfrak{w}}}^{(n-1)} - \widetilde{\boldsymbol{\mathfrak{w}}}^{(n+1)}) \rangle   + \frac{4(\Phi_{\max}(\boldsymbol{L}_-)^2 + \Phi_{\max}(\boldsymbol{L}_+)^2)}{\Phi_{\min}(\boldsymbol{L}_-)} ||\boldsymbol{\xi}^{(n+1)}||^2_2 \notag\\
    &\quad + \langle \widetilde{\boldsymbol{\mathfrak{w}}}^{(n)} - \boldsymbol{\mathfrak{w}}^{*}, \frac{2}{\rho} \boldsymbol{D}^{(n+1)} \otimes \boldsymbol{I}_P (\widetilde{\boldsymbol{\mathfrak{w}}}^{(n)} - \widetilde{\boldsymbol{\mathfrak{w}}}^{(n+1)})\rangle \notag
\end{align}
where $\boldsymbol{q}^*=[\boldsymbol{r}^T,(\boldsymbol{\mathfrak{w}}^*)^T]$.}
\begin{proof}
See Appendix D.
\end{proof}
Following the result of Lemma II, we can establish the following theorem from which we will derive the converge results. 

\noindent \textbf{Theorem II.} \textit{Under Assumption 2, and given the final iteration $T > 0$, we can bound the expected error of the zCDP-NFL algorithm as
\begin{align}
    &\mathbb{E}[f(\hat{\boldsymbol{\mathfrak{w}}}^{(T)}) - f(\boldsymbol{\mathfrak{w}}^*)]  \leqslant \frac{\rho}{T} \sum_{n=1}^{T} \Bigl( - 2 \langle \boldsymbol{Q} \widetilde{\boldsymbol{\mathfrak{w}}}^{(n)}, \boldsymbol{Q} \widetilde{\boldsymbol{\mathfrak{w}}}^{(n+1)} \rangle - ||\widetilde{\boldsymbol{\mathfrak{w}}}^{(n)} - \widetilde{\boldsymbol{\mathfrak{w}}}^{(n-1)}||_{\frac{\boldsymbol{L}_+}{2}}^2 \notag\\
    &\qquad + \langle \widetilde{\boldsymbol{\mathfrak{w}}}^{(n)} - \boldsymbol{\mathfrak{w}}^{*}, \frac{2}{\rho} \boldsymbol{D}^{(n+1)} \otimes \boldsymbol{I}_P (\widetilde{\boldsymbol{\mathfrak{w}}}^{(n)} - \widetilde{\boldsymbol{\mathfrak{w}}}^{(n+1)})\rangle \notag\\
    &\qquad - \langle\boldsymbol{\mathfrak{w}}^{*}, \boldsymbol{L}_+(2 \widetilde{\boldsymbol{\mathfrak{w}}}^{(n)} - \widetilde{\boldsymbol{\mathfrak{w}}}^{(n-1)} - \widetilde{\boldsymbol{\mathfrak{w}}}^{(n+1)}) \rangle 
     +  ||\widetilde{\boldsymbol{\mathfrak{w}}}^{(n+1)} - \widetilde{\boldsymbol{\mathfrak{w}}}^{(n)}||_{\boldsymbol{L}_+}^2 \Bigr) \notag\\
    &\quad + \frac{1}{T} \frac{\rho P 4(\Phi_{\max}(\boldsymbol{L}_-)^2 + \Phi_{\max}(\boldsymbol{L}_+)^2) \sum_{k=1}^K \sigma^{2 (0)}_k} {\Phi_{\min}(\boldsymbol{L}_-) (1 - \tau)} \notag\\
    &\; + \frac{\langle \widetilde{\boldsymbol{\mathfrak{w}}}^{(1)}, \boldsymbol{L}_+(\widetilde{\boldsymbol{\mathfrak{w}}}^{(1)} - \widetilde{\boldsymbol{\mathfrak{w}}}^{(0)})\rangle}{T}+ \frac{\rho || \boldsymbol{Q} \widetilde{\boldsymbol{\mathfrak{w}}}^{(0)}||_2^2}{T} + \frac{ \rho ||\widetilde{\boldsymbol{\mathfrak{w}}}^{(0)} - \boldsymbol{\mathfrak{w}}^*||_{\frac{ \boldsymbol{L}_-}{2}}^2}{T}.     \label{Theorem.II}
\end{align}
where $\hat{\boldsymbol{\mathfrak{w}}}^{(T)} = \frac{1}{T} \sum_{n=1}^{T} \widetilde{\boldsymbol{\mathfrak{w}}}^{(n)}$, and the expectation is taken with respect to the noise. Since $\boldsymbol{\mathfrak{w}}^*$ is the optimal solution, $\mathbb{E}[f(\hat{\boldsymbol{\mathfrak{w}}}^{(T)}) - f(\boldsymbol{\mathfrak{w}}^*)]$ is positive.}
\begin{proof}
See Appendix E.
\end{proof}

\subsection{Convergence Properties}

We can derive three important results from Theorem II. The first is that the zCDP-NFL algorithm converges to the exact solution of \eqref{distribERM}. The second is the rate of this convergence. The third result is the privacy accuracy trade-off bound of the algorithm. First, we define the required assumptions for convergence.

\noindent \textbf{Assumption 3.} \textit{We require that $\Lim{n \to +\infty} \eta_k^{(n)} = 0, \forall k \in \mathcal{C}$. This will enforce the asymptotic stability of the local estimates.}

\noindent \textbf{Theorem III.} \textit{Under Assumptions 2 and 3, the zCDP-NFL algorithm defined by the steps \eqref{PrimalUpdate}-\eqref{DualUpdate}, converges to the exact solution.}
\begin{proof}
We can simplify the result of Theorem II into the following:
\begin{align}
    &\mathbb{E}[f(\hat{\boldsymbol{\mathfrak{w}}}^{(T)}) - f(\boldsymbol{\mathfrak{w}}^*)]   \leqslant \frac{\rho}{T} \sum_{n=1}^{T} \Bigl( - 2 \langle \boldsymbol{Q} \widetilde{\boldsymbol{\mathfrak{w}}}^{(n)}, \boldsymbol{Q} \widetilde{\boldsymbol{\mathfrak{w}}}^{(n+1)} \rangle \notag\\
    &\qquad + \langle \widetilde{\boldsymbol{\mathfrak{w}}}^{(n)} - \boldsymbol{\mathfrak{w}}^{*}, \frac{2}{\rho} \boldsymbol{D}^{(n+1)} \otimes \boldsymbol{I}_P (\widetilde{\boldsymbol{\mathfrak{w}}}^{(n)} - \widetilde{\boldsymbol{\mathfrak{w}}}^{(n+1)})\rangle   + ||\widetilde{\boldsymbol{\mathfrak{w}}}^{(n+1)} - \widetilde{\boldsymbol{\mathfrak{w}}}^{(n)}||_{\boldsymbol{L}_+}^2 \Bigr) \notag\\
    &\quad + \frac{1}{T} \frac{\rho P 4(\Phi_{\max}(\boldsymbol{L}_-)^2 + \Phi_{\max}(\boldsymbol{L}_+)^2) \sum_{k=1}^K \sigma^{2 (0)}_k} {\Phi_{\min}(\boldsymbol{L}_-) (1 - \tau)} \notag\\
    &\, + \frac{\langle \widetilde{\boldsymbol{\mathfrak{w}}}^{(1)}, \boldsymbol{L}_+(\widetilde{\boldsymbol{\mathfrak{w}}}^{(1)} - \widetilde{\boldsymbol{\mathfrak{w}}}^{(0)})\rangle}{T} + \frac{\rho || \boldsymbol{Q} \widetilde{\boldsymbol{\mathfrak{w}}}^{(0)}||_2^2}{T} + \frac{ \rho ||\widetilde{\boldsymbol{\mathfrak{w}}}^{(0)} - \boldsymbol{\mathfrak{w}}^*||_{\frac{ \boldsymbol{L}_-}{2}}^2}{T}.     \label{Th3}
\end{align}
We will consider the terms separately in their order of appearance. We first prove that $ \Lim{n \to +\infty} \frac{\rho \sum_{n=1}^{T} - 2 \langle \boldsymbol{Q} \widetilde{\boldsymbol{\mathfrak{w}}}^{(n)}, \boldsymbol{Q} \widetilde{\boldsymbol{\mathfrak{w}}}^{(n+1)} \rangle}{T} = 0$. 

We can now note that 
\begin{align}
    & - 2 \langle \boldsymbol{Q} \widetilde{\boldsymbol{\mathfrak{w}}}^{(n)}, \boldsymbol{Q} \widetilde{\boldsymbol{\mathfrak{w}}}^{(n+1)} \rangle  
     = - \langle \boldsymbol{Q} \widetilde{\boldsymbol{\mathfrak{w}}}^{(n)}, \boldsymbol{Q} \widetilde{\boldsymbol{\mathfrak{w}}}^{(n)} \rangle - \langle \boldsymbol{Q} \widetilde{\boldsymbol{\mathfrak{w}}}^{(n)}, \boldsymbol{Q} (\widetilde{\boldsymbol{\mathfrak{w}}}^{(n+1)} - \widetilde{\boldsymbol{\mathfrak{w}}}^{(n)}) \rangle \notag\\
    &\quad - \langle \boldsymbol{Q} \widetilde{\boldsymbol{\mathfrak{w}}}^{(n+1)}, \boldsymbol{Q} \widetilde{\boldsymbol{\mathfrak{w}}}^{(n+1)} \rangle - \langle \boldsymbol{Q} (\widetilde{\boldsymbol{\mathfrak{w}}}^{(n)} - \widetilde{\boldsymbol{\mathfrak{w}}}^{(n+1)}), \boldsymbol{Q} \widetilde{\boldsymbol{\mathfrak{w}}}^{(n+1)} \rangle \notag\\
    &= - ||\boldsymbol{Q} \widetilde{\boldsymbol{\mathfrak{w}}}^{(n)}||^2_2 - ||\boldsymbol{Q} \widetilde{\boldsymbol{\mathfrak{w}}}^{(n+1)}||^2_2 + ||\boldsymbol{Q} (\widetilde{\boldsymbol{\mathfrak{w}}}^{(n+1)} - \widetilde{\boldsymbol{\mathfrak{w}}}^{(n)})||^2_2 \notag\\
    &\leqslant - ||\boldsymbol{Q} \widetilde{\boldsymbol{\mathfrak{w}}}^{(n)}||^2_2 - ||\boldsymbol{Q} \widetilde{\boldsymbol{\mathfrak{w}}}^{(n+1)}||^2_2 + ||\boldsymbol{Q}||_2^2 ||\widetilde{\boldsymbol{\mathfrak{w}}}^{(n+1)} - \widetilde{\boldsymbol{\mathfrak{w}}}^{(n)}||^2_2.  
    \label{PbAnalysis}
\end{align}

As seen in \eqref{PrimalUpdate}, $ \boldsymbol{\mathfrak{w}}^{(n+1)} $ minimizes a function where all terms are bounded except the term $\frac{\big\|\boldsymbol{\mathfrak{w}}-\widetilde{\boldsymbol{\mathfrak{w}}}^{(n)}\big\|^2}{2\eta_k^{(n+1)}}$. Therefore, under Assumption 3, $\Lim{n \to +\infty} ||\boldsymbol{\mathfrak{w}}^{(n+1)} - \widetilde{\boldsymbol{\mathfrak{w}}}^{(n)}||^2_2 = 0$. Since $\widetilde{\boldsymbol{\mathfrak{w}}}^{(n+1)}$ is defined as $\widetilde{\boldsymbol{\mathfrak{w}}}^{(n+1)} = \boldsymbol{\mathfrak{w}}^{(n+1)} + \boldsymbol{\xi}^{(n+1)}$ with $\Lim{n \to +\infty} ||\boldsymbol{\xi}^{(n+1)}|| = 0$, we have $\Lim{n \to +\infty} ||\widetilde{\boldsymbol{\mathfrak{w}}}^{(n+1)} - \widetilde{\boldsymbol{\mathfrak{w}}}^{(n)}||^2_2 = 0$. \\

This implies that $- 2 \langle \boldsymbol{Q} \widetilde{\boldsymbol{\mathfrak{w}}}^{(n)}, \boldsymbol{Q} \widetilde{\boldsymbol{\mathfrak{w}}}^{(n+1)} \rangle$ is bounded by a series converging to $0$. Therefore, since $\mathbb{E}[\hat{f}(\hat{\boldsymbol{\mathfrak{w}}}^{(T)}) - \hat{f}(\boldsymbol{\mathfrak{w}}^*)]$ is positive, $ \frac{\rho \sum_{n=1}^{T} - 2 \langle \boldsymbol{Q} \widetilde{\boldsymbol{\mathfrak{w}}}^{(n)}, \boldsymbol{Q} \widetilde{\boldsymbol{\mathfrak{w}}}^{(n+1)} \rangle}{T}$ converges to $0$. 

Next, under Assumption 3, we have $ \Lim{n \to +\infty} \frac{\rho \sum_{n=1}^{T} \langle \widetilde{\boldsymbol{\mathfrak{w}}}^{(n)} - \boldsymbol{\mathfrak{w}}^{*}, \frac{2}{\rho} \boldsymbol{D}^{(n+1)} \otimes \boldsymbol{I}_P (\widetilde{\boldsymbol{\mathfrak{w}}}^{(n)} - \widetilde{\boldsymbol{\mathfrak{w}}}^{(n+1)})\rangle}{T} = 0 $ since $[\boldsymbol{D}^{(n+1)}]_{k,k} = \frac{1}{\sqrt{2 \eta_k^{(n+1)}}}$. We now consider $\frac{\rho}{T} \sum_{n=1}^{T} ||\widetilde{\boldsymbol{\mathfrak{w}}}^{(n+1)} - \widetilde{\boldsymbol{\mathfrak{w}}}^{(n)}||_{\boldsymbol{L}_+}^2$. As we have shown that $\Lim{n \to +\infty} ||\widetilde{\boldsymbol{\mathfrak{w}}}^{(n+1)} - \widetilde{\boldsymbol{\mathfrak{w}}}^{(n)}||^2_2= 0$, we have $\Lim{n \to +\infty} ||\widetilde{\boldsymbol{\mathfrak{w}}}^{(n+1)} - \widetilde{\boldsymbol{\mathfrak{w}}}^{(n)}||_{\boldsymbol{L}_+}^2  = 0$, and therefore, the sum $\sum_{n=1}^{T} ||\widetilde{\boldsymbol{\mathfrak{w}}}^{(n+1)} - \widetilde{\boldsymbol{\mathfrak{w}}}^{(n)}||_{\boldsymbol{L}_+}^2$ is a Cauchy sequence. Hence we have $\Lim{n \to +\infty} \frac{\rho}{T} \sum_{n=1}^{T} ||\widetilde{\boldsymbol{\mathfrak{w}}}^{(n+1)} - \widetilde{\boldsymbol{\mathfrak{w}}}^{(n)}||_{\boldsymbol{L}_+}^2 = 0$. 


Finally, the terms outside of the summation trivially converge to $0$ as $T \to +\infty$. This concludes the proof.
\end{proof}

We now introduce the required assumption to establish the convergence rate of the algorithm.

\noindent \textbf{Assumption 4.} \textit{The $\eta_k^{(n)}, k \in \mathcal{C}$ are chosen such that $||\boldsymbol{D}^{(n+1)}||^2_2$ is a convergent series. This assumption, stronger than Assumption 3, is necessary to guarantee the exponential stability of the local estimates.}


\noindent \textbf{Theorem IV.} \textit{Under Assumptions 2 and 4, the zCDP-NFL algorithm converges with a rate of $\mathcal{O}(1/n)$ iterations.}

\begin{proof}
In the following, we assume that the optimal solution $\boldsymbol{\mathfrak{w}}^{*} \ne \boldsymbol{0}$. If $\boldsymbol{\mathfrak{w}}^{*} = \boldsymbol{0}$, one could add a nonzero artificial dimension and proceed.

In order to prove this result, we will show that the expectation of the error is bounded by a bounded term divided by $T$. Notably, we will show that the sum in \eqref{Th3} converges. We consider the terms in their order of appearance in Theorem II. We will also use the result of Theorem III, $\Lim{n \to +\infty} \widetilde{\boldsymbol{\mathfrak{w}}}^{(n)} = \boldsymbol{\mathfrak{w}}^{*}$, for which Assumption 3 is satisfied by Assumption 4.

To begin, we consider $\frac{\rho}{T} \sum_{n=1}^{T} \Bigl( - 2 \langle \boldsymbol{Q} \widetilde{\boldsymbol{\mathfrak{w}}}^{(n)}, \boldsymbol{Q} \widetilde{\boldsymbol{\mathfrak{w}}}^{(n+1)} \rangle - ||\widetilde{\boldsymbol{\mathfrak{w}}}^{(n)} - \widetilde{\boldsymbol{\mathfrak{w}}}^{(n-1)}||_{\frac{\boldsymbol{L}_+}{2}}^2 \Bigr)$. Since $\widetilde{\boldsymbol{\mathfrak{w}}}^{(n)}$ converges to $\widetilde{\boldsymbol{\mathfrak{w}}}^{*}$, $- 2 \langle \boldsymbol{Q} \widetilde{\boldsymbol{\mathfrak{w}}}^{(n)}, \boldsymbol{Q} \widetilde{\boldsymbol{\mathfrak{w}}}^{(n+1)} \rangle$ converges to $ - 2 ||\boldsymbol{Q} \boldsymbol{\mathfrak{w}}^{*}||_2^2 $ that is strictly negative. Therefore, there exist an iteration $n_0$ after which all terms $- 2 \langle \boldsymbol{Q} \widetilde{\boldsymbol{\mathfrak{w}}}^{(n)}, \boldsymbol{Q} \widetilde{\boldsymbol{\mathfrak{w}}}^{(n+1)} \rangle$ are negative. Hence, $\frac{\rho}{T} \sum_{n=1}^{T} \Bigl( - 2 \langle \boldsymbol{Q} \widetilde{\boldsymbol{\mathfrak{w}}}^{(n)},$ $ \boldsymbol{Q} \widetilde{\boldsymbol{\mathfrak{w}}}^{(n+1)} \rangle - ||\widetilde{\boldsymbol{\mathfrak{w}}}^{(n)} - \widetilde{\boldsymbol{\mathfrak{w}}}^{(n-1)}||_{\frac{\boldsymbol{L}_+}{2}}^2 \Bigr) \leqslant \frac{\rho}{T} \sum_{n=1}^{n_0} - 2 \langle \boldsymbol{Q} \widetilde{\boldsymbol{\mathfrak{w}}}^{(n)}, \boldsymbol{Q} \widetilde{\boldsymbol{\mathfrak{w}}}^{(n+1)} \rangle$.

Next, we can bound $\langle \widetilde{\boldsymbol{\mathfrak{w}}}^{(n)} - \boldsymbol{\mathfrak{w}}^{*}, \frac{2}{\rho} \boldsymbol{D}^{(n+1)} \otimes \boldsymbol{I}_P (\widetilde{\boldsymbol{\mathfrak{w}}}^{(n)} - \widetilde{\boldsymbol{\mathfrak{w}}}^{(n+1)})\rangle$ by $||\widetilde{\boldsymbol{\mathfrak{w}}}^{(n)} - \boldsymbol{\mathfrak{w}}^{*}||^2_2 \frac{2}{\rho} ||\boldsymbol{D}^{(n+1)}||^2_2 ||\boldsymbol{I}_P||^2_2 ||\widetilde{\boldsymbol{\mathfrak{w}}}^{(n)} - \widetilde{\boldsymbol{\mathfrak{w}}}^{(n+1)}||^2_2$ and thus $\rho \sum_{n=1}^{T} \langle \widetilde{\boldsymbol{\mathfrak{w}}}^{(n)} - \boldsymbol{\mathfrak{w}}^{*}, \frac{2}{\rho} \boldsymbol{D}^{(n+1)} \otimes \boldsymbol{I}_P (\widetilde{\boldsymbol{\mathfrak{w}}}^{(n)} - \widetilde{\boldsymbol{\mathfrak{w}}}^{(n+1)})\rangle$ by $2 \sum_{n=1}^{T} ||\widetilde{\boldsymbol{\mathfrak{w}}}^{(n)} - \boldsymbol{\mathfrak{w}}^{*}||^2_2||\boldsymbol{D}^{(n+1)}||^2_2 ||\widetilde{\boldsymbol{\mathfrak{w}}}^{(n)} - \widetilde{\boldsymbol{\mathfrak{w}}}^{(n+1)}||^2_2$. Using $\Lim{n \to +\infty} \widetilde{\boldsymbol{\mathfrak{w}}}^{(n)} = \boldsymbol{\mathfrak{w}}^{*}$, there exist constants $\alpha_0$ and $\alpha_1$ such that $\forall n \in \mathbb{N}, ||\widetilde{\boldsymbol{\mathfrak{w}}}^{(n)} - \boldsymbol{\mathfrak{w}}^{*}||^2_2 \leqslant \alpha_0$ and $||\widetilde{\boldsymbol{\mathfrak{w}}}^{(n)} - \widetilde{\boldsymbol{\mathfrak{w}}}^{(n+1)}||^2_2 \leqslant \alpha_1$. This leads to $2 \sum_{n=1}^{T} ||\widetilde{\boldsymbol{\mathfrak{w}}}^{(n)} - \boldsymbol{\mathfrak{w}}^{*}||^2_2||\boldsymbol{D}^{(n+1)}||^2_2 ||\widetilde{\boldsymbol{\mathfrak{w}}}^{(n)} - \widetilde{\boldsymbol{\mathfrak{w}}}^{(n+1)}||^2_2 \leqslant 2 \alpha_0 \alpha_1 \sum_{n=1}^{T} ||\boldsymbol{D}^{(n+1)}||^2_2$, which is a convergent series under Assumption 4. Therefore, $\frac{\rho}{T} \sum_{n=1}^{T} \langle \widetilde{\boldsymbol{\mathfrak{w}}}^{(n)} - \boldsymbol{\mathfrak{w}}^{*}, \frac{2}{\rho} \boldsymbol{D}^{(n+1)} \otimes \boldsymbol{I}_P (\widetilde{\boldsymbol{\mathfrak{w}}}^{(n)} - \widetilde{\boldsymbol{\mathfrak{w}}}^{(n+1)})\rangle$ converges to zero in $O(1/n)$.

We can bound $\langle\boldsymbol{\mathfrak{w}}^{*}, \boldsymbol{L}_+(2 \widetilde{\boldsymbol{\mathfrak{w}}}^{(n)} - \widetilde{\boldsymbol{\mathfrak{w}}}^{(n-1)} - \widetilde{\boldsymbol{\mathfrak{w}}}^{(n+1)}) \rangle$ by $2 ||\boldsymbol{\mathfrak{w}}^{*}||_2^2 ||\boldsymbol{L}_+||_2^2 (|| \widetilde{\boldsymbol{\mathfrak{w}}}^{(n+1)} - \widetilde{\boldsymbol{\mathfrak{w}}}^{(n)}||_2^2 + || \widetilde{\boldsymbol{\mathfrak{w}}}^{(n)} - \widetilde{\boldsymbol{\mathfrak{w}}}^{(n-1)} ||_2^2) $. Using $\Lim{n \to +\infty} \widetilde{\boldsymbol{\mathfrak{w}}}^{(n)} = \boldsymbol{\mathfrak{w}}^{*}$, the series $\sum_{n=1}^{\infty} ||\widetilde{\boldsymbol{\mathfrak{w}}}^{(n+1)} - \widetilde{\boldsymbol{\mathfrak{w}}}^{(n)}||_2^2$ and $\sum_{n=1}^{\infty} ||\widetilde{\boldsymbol{\mathfrak{w}}}^{(n+1)} - \widetilde{\boldsymbol{\mathfrak{w}}}^{(n)}||_{\boldsymbol{L}_+}^2$ converge to values that we denote $\alpha_2$ and $\alpha_3$, respectively. We have $ \sum_{n=1}^{T} ( - \langle\boldsymbol{\mathfrak{w}}^{*}, \boldsymbol{L}_+(2 \widetilde{\boldsymbol{\mathfrak{w}}}^{(n)} - \widetilde{\boldsymbol{\mathfrak{w}}}^{(n-1)} - \widetilde{\boldsymbol{\mathfrak{w}}}^{(n+1)}) \rangle + ||\widetilde{\boldsymbol{\mathfrak{w}}}^{(n+1)} - \widetilde{\boldsymbol{\mathfrak{w}}}^{(n)}||_{\boldsymbol{L}_+}^2) \leqslant 4 ||\boldsymbol{\mathfrak{w}}^{*}||_2^2 ||\boldsymbol{L}_+||_2^2 \alpha_2 + \alpha_3$.

Finally, we prove that all terms outside of the sum are bounded by a constant with respect to $T$. The only one requiring further analysis is $\langle \widetilde{\boldsymbol{\mathfrak{w}}}^{(1)}, \boldsymbol{L}_+(\widetilde{\boldsymbol{\mathfrak{w}}}^{(1)} - \widetilde{\boldsymbol{\mathfrak{w}}}^{(0)})\rangle$ and it can be bounded by $||\widetilde{\boldsymbol{\mathfrak{w}}}^{(1)}||^2_2 ||\boldsymbol{L}_+( \widetilde{\boldsymbol{\mathfrak{w}}}^{(1)} - \widetilde{\boldsymbol{\mathfrak{w}}}^{(0)})||^2_2$.

Each term has either been bounded by a constant with respect to $T$, divided by $T$; this concludes the proof.
\end{proof}

\noindent \textbf{Remark:} In practice, Assumption 4 can be relaxed in most cases.

\subsection{Privacy Accuracy Trade-off}

The last result established by Theorem II is the privacy accuracy trade-off bound. The privacy accuracy trade-off quantifies how ensuring more privacy deteriorates the accuracy of the algorithm and is one of the most important parameters of a privacy-preserving algorithm. Under Assumption 4, we can reformulate \eqref{Theorem.II} as
\begin{align}
    &\mathbb{E}[f(\hat{\boldsymbol{\mathfrak{w}}}^{(T)}) - f(\boldsymbol{\mathfrak{w}}^*)] \leqslant \frac{\alpha}{T} + \frac{\alpha_{\boldsymbol{\xi}}}{T} \frac{\sum_{k=1}^K \sigma^{2}_k(0)} {1 - \tau}
\end{align}
where $\alpha$ is a constant with respect to $T$ and the noise perturbation and $\alpha_{\boldsymbol{\xi}} = \frac{\rho P 4(\Phi_{\max}(\boldsymbol{L}_-)^2 + \Phi_{\max}(\boldsymbol{L}_+)^2)} {\Phi_{\min}(\boldsymbol{L}_-)}$.

By combining this result with Theorem I, we obtain
\begin{align}
    &\mathbb{E}[f(\hat{\boldsymbol{\mathfrak{w}}}^{(T)}) - f(\boldsymbol{\mathfrak{w}}^*)] \leqslant \frac{\alpha}{T} + \frac{\alpha_{\boldsymbol{\xi}}}{T} \frac{\sum_{k=1}^K \frac{\Delta_{k,2}^2(0)}{2\varphi_{k}^{(1)}}}{1-\tau}
\end{align}

In the common case where the privacy parameter $\varphi_{k}^{(1)}$ is identical for all clients, i.e., $\varphi_{k}^{(1)} = \varphi^{(1)}, \forall k \in \mathcal{C}$, we have
\begin{align}
    &\mathbb{E}[f(\hat{\boldsymbol{\mathfrak{w}}}^{(T)}) - f(\boldsymbol{\mathfrak{w}}^*)] \leqslant \frac{\alpha}{T} + \frac{\alpha_{\boldsymbol{\xi}}}{T} \frac{\sum_{k=1}^K \Delta_{k,2}^2(0)}{2 K \varphi^{(1)} (1-\tau)}
\end{align}
With this result, we see that ensuring more privacy, which can be done by decreasing $\varphi^{(1)}$ or having $\tau$ closer to $1$, would result in a less restrictive convergence bound for the algorithm.

\section{Numerical Simulations}

This section presents simulation results to evaluate the performance and privacy accuracy trade-off of the proposed zCDP-NFL. To compare the different DP implementations, we introduce the $(\epsilon, \delta)$DP-NFL, identical to zCDP-NFL except for the fact that it uses conventional $(\epsilon, \delta)$-DP rules to control the local noise perturbation. On nonsmooth objective functions, we benchmark the proposed algorithm against conventional subgradient-based networked FL, as presented in \cite{Nedic2009}, modified to use zCDP and denoted zCDP-grad-NFL. On smooth objective functions, we benchmark the proposed algorithm against the existing distributed learning algorithm P-ADMM, presented in \cite{Zhuhan} that uses a networked FL architecture, zCDP, and the ADMM, but is constrained to smooth objective functions. In the following, we consider the elastic net, least absolute deviation, and ridge regression problems, all presented in \cite{Boyd2010}.

For a fair comparison, the algorithms are tuned to provide the same total privacy guarantees throughout the computation - this is made possible by the corollary of Theorem I. This corollary provides $(\epsilon, \delta)$-DP guarantees for an algorithm using zCDP with both $\varphi^{(1)}$ and $\tau$ as parameters.  Furthermore, the algorithms are tuned to observe identical initial convergence speeds when possible.

In the following, we consider a network with a random topology comprising $K = 50$ nodes, where each node connects to $3$ other nodes on average. Each node $k$ possesses $M_k=50$ local noisy observations of the unknown parameter ${\bf w}$ of dimension $P=8$. The proposed simulations are performed on synthetic data and performance is evaluated by computing the normalized error defined as $\sum_{k=1}^{K}{||{\bf w}_{k}^{(n)}-{\bf w}_c||^2}/{||{\bf w}_c||^2}$, ${\bf w}_c$ being the centralized solution obtained by the CVX toolbox \cite{cvx}.

\begin{figure*}[t]
    \centering
    \includegraphics[width=0.45\textwidth]{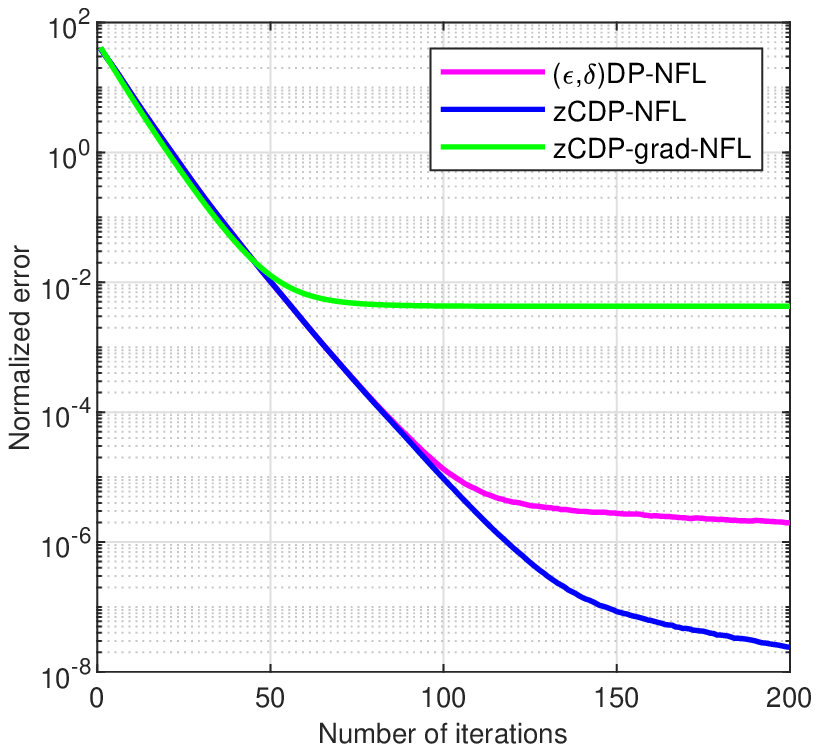}
    \includegraphics[width=0.45\textwidth]{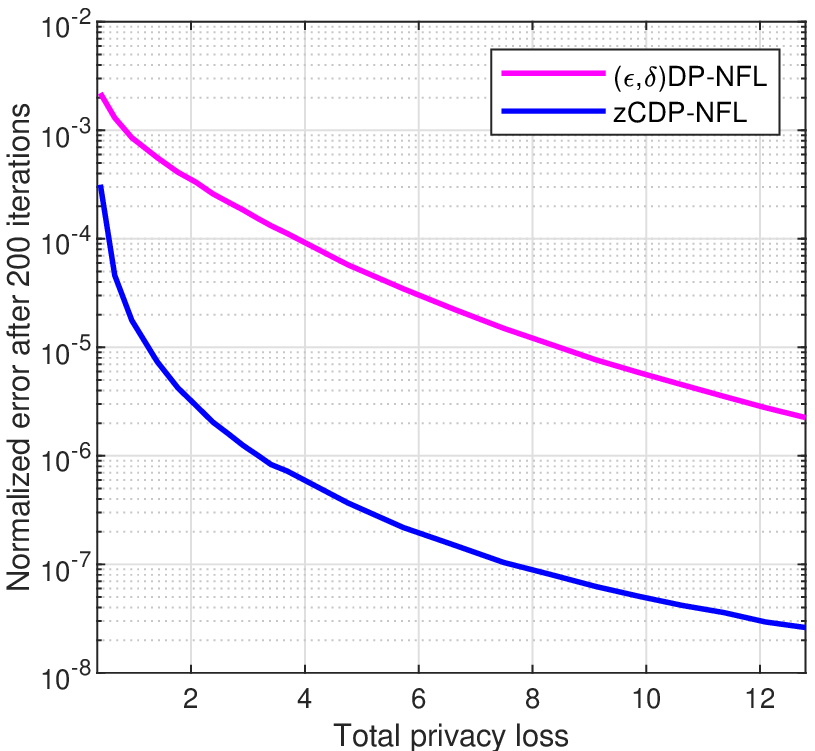}
    \caption{Learning curves (left) and privacy-accuracy trade-off (right) on the elastic net problem.}
    \label{Convergence_EN}
\end{figure*}

Figure \ref{Convergence_EN} contains the learning curve, i.e., normalized error versus iteration index, and privacy-accuracy trade-off, i.e., normalized error after $200$ iterations versus total privacy loss, on the elastic net problem, defined by a smooth loss function $\ell(\mathbf{X}_k,\mathbf{y}_k;{\bf w}_k)=||\mathbf{X}_k{\bf w}_k-\mathbf{y}_k||^2$ and a nonsmooth regularizer function $R({\bf w}_k)= \lambda_1 ||{\bf w}_k||_1 + \lambda_2 ||{\bf w}_k||^2$ with $\lambda_1 = 0.001||\mathbf{X}^\mathsf{T}\mathbf{y}||_{\infty}$, as in \cite{Boyd2010}, and $\lambda_2 = 1$. We observe that the proposed ADMM-based algorithm significantly outperforms the subgradient-based algorithm. Furthermore, we can see that the use of the zCDP notion as opposed to the $(\epsilon, \delta)$-DP notion allows for better accuracy given the same privacy budget throughout the computation. This same fact is illustrated for various total privacy losses on the privacy-accuracy trade-off curve.



\begin{figure*}[t]
    \centering
    \includegraphics[width=0.45\textwidth]{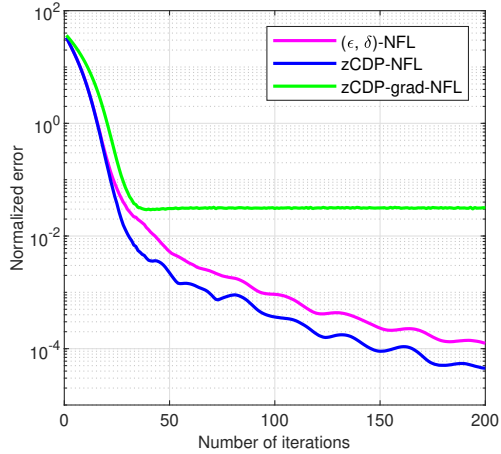}
    \caption{Learning curve on the least absolute deviation problem.}
    \label{Convergence_LAD}
\end{figure*}

Figure \ref{Convergence_LAD} contains the learning curve on the least absolute deviation problem, solely composed of a nonsmooth loss $\ell(\mathbf{X}_k,\mathbf{y}_k;{\bf w}_k) = ||\mathbf{X}_k{\bf w}_k-\mathbf{y}_k||_1$. We observe once again that the proposed algorithm significantly outperforms its subgradient-based counterpart, as the learning rate required for the zCDP-grad-NFL algorithm to attain a similar initial convergence speed to the ADMM-based algorithms does not allow it to reach high accuracy. In addition, we observe that the use of zCDP allows for better accuracy than $(\epsilon, \delta)$-DP. In the end, the proposed zCDP-NFL algorithm significantly improves over existing methods on nonsmooth objectives.

\begin{figure*}[t]
    \centering
    \includegraphics[width=0.45\textwidth]{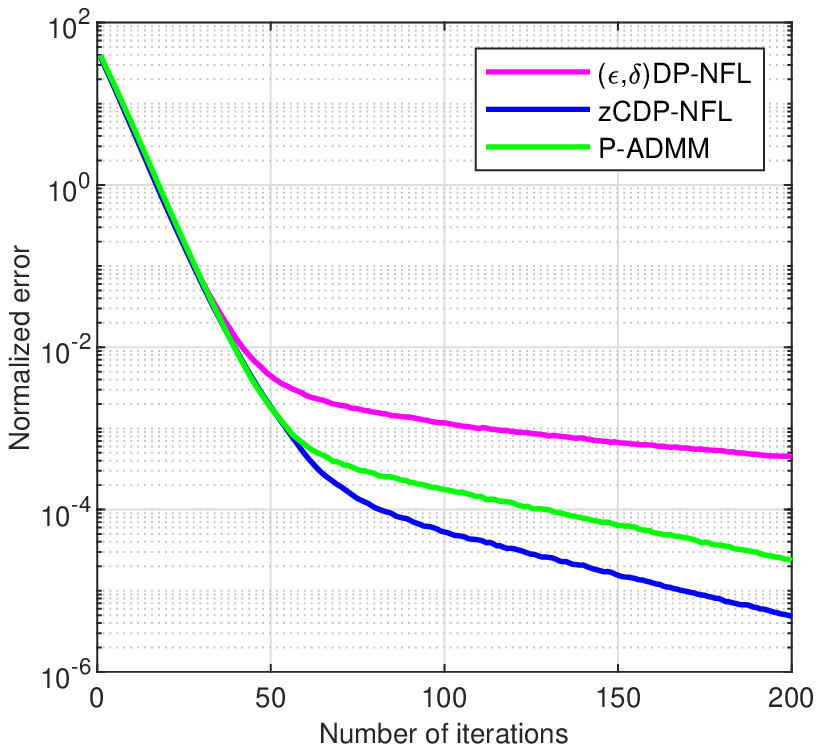}
    \includegraphics[width=0.45\textwidth]{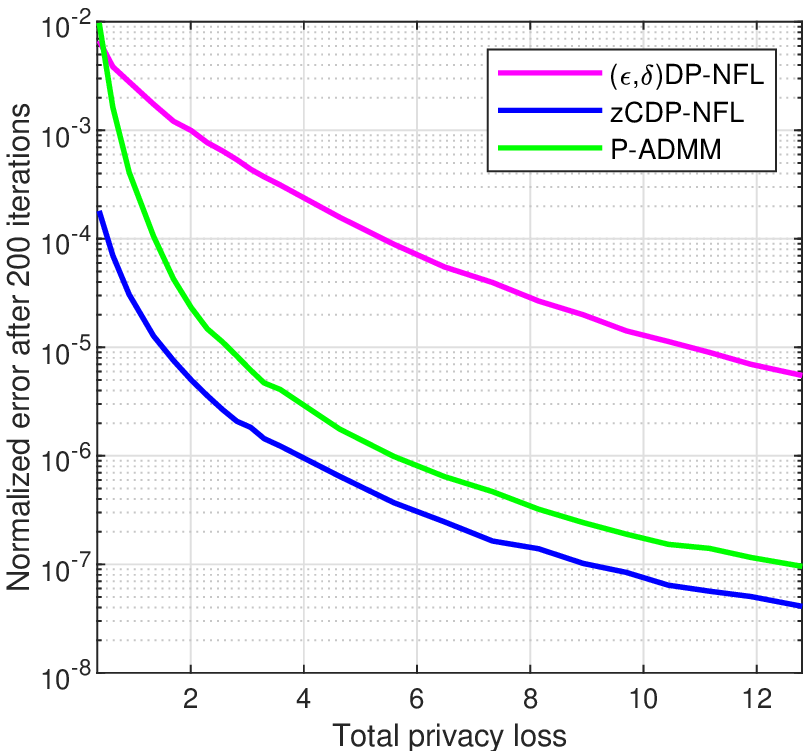}
    \caption{Learning curves (left) and privacy-accuracy trade-off (right) on the smooth ridge regression problem.}
    \label{Convergence_RR}
\end{figure*}

Figure \ref{Convergence_RR} contains the learning curve and privacy-accuracy trade-off on the ridge regression problem, defined by a smooth loss function $\ell(\mathbf{X}_k,\mathbf{y}_k;{\bf w}_k) = ||\mathbf{X}_k{\bf w}_k-\mathbf{y}_k||^2$ and a smooth regularizer function $R({\bf w}_k)=||{\bf w}_k||^2$. Since this objective is smooth, we can present the performance of the P-ADMM algorithm. We observe that the proposed zCDP-NFL algorithm slightly outperforms the P-ADMM algorithm on smooth objectives, despite the inexact ADMM update required to handle nonsmooth objectives. As can be observed on the privacy-accuracy trade-off curve, this is the case for all total privacy losses. This better performance is due to the use of the time-varying step size $\eta$ in the proposed algorithm. Even though the zCPD-NFL algorithm is designed for nonsmooth objective functions, it offers very good performances on smooth objective functions.

\section{Conclusions and Future Directions}
The proposed zCDP-NFL algorithm is a networked, privacy-preserving algorithm that accommodates nonsmooth and non-strongly convex objective functions. Each client is protected by local differential privacy with guarantees provided in terms of zCDP. We provided mathematical proofs of the privacy guarantee and convergence to the optimal point in $O(1/n)$, as well as an analysis of the privacy-accuracy trade-off to quantify the accuracy loss caused by increased privacy. Numerical simulations show that the proposed zCDP-NFL algorithm significantly outperforms existing networked algorithms on nonsmooth objectives while offering very good performances on smooth objectives. Future work includes communication efficient implementations and robustness to model poisoning.


\appendix

\section{Proof of Lemma I}
\label{sec:sample:appendix}

\begin{proof}

We consider two neighboring data sets $\mathcal{D}$ and $\mathcal{D}'$ and their respective primal updates for the client $k$ whose data set contains the difference. We will denote $\mathcal{D}_k$ and $\mathcal{D}'_k$ the local data set of client $k$ corresponding to the use of $\mathcal{D}$ and $\mathcal{D}'$, respectively. Moreover, we will denote ${\bf w}_{k,\mathcal{D}_k}^{(n)}$ and ${\bf w}_{k,\mathcal{D}_k'}^{(n)}$ the estimates computed by client $k$ using the data sets $\mathcal{D}_k$ and $\mathcal{D}'_k$, respectively. These can be computed as
\begin{align}
    &{\bf w}_{k,\mathcal{D}_k}^{(n)}=  \frac{1}{2 \rho |\mathcal{N}_k| + \frac{1}{\eta^{(n)}}} \Bigr( \frac{\widetilde{{\bf w}}_{k}^{(n-1)}}{\eta^{(n)}} + \frac{\rho}{2} \sum_{i \in \mathcal{N}_k} (\widetilde{{\bf w}}_{k}^{(n-1)} - \widetilde{{\bf w}}_{i}^{(n-1)}) \notag\\
    & \quad \quad + \frac{\boldsymbol{\gamma}_k^{(n-1)}}{2}   -\sum_{j=1}^{M_k} \frac{\ell'(\mathbf{x}_{k,j},y_{k,j}; \widetilde{{\bf w}}_k)}{M_k} - \frac{\lambda R'({\bf w}_k)}{K} \Bigr), 
\end{align}
\begin{align}
    &{\bf w}_{k,\mathcal{D}_k'}^{(n)}= \frac{1}{2 \rho |\mathcal{N}_k| + \frac{1}{\eta^{(n)}}} \Bigr( \frac{\widetilde{{\bf w}}_{k}^{(n-1)}}{\eta^{(n)}} + \frac{\rho}{2} \sum_{i \in \mathcal{N}_k} (\widetilde{{\bf w}}_{k}^{(n-1)} - \widetilde{{\bf w}}_{i}^{(n-1)}) + \frac{\boldsymbol{\gamma}_k^{(n-1)}}{2}  \notag\\
    & \quad \quad -\sum_{j=1}^{M_k-1} \frac{\ell'(\mathbf{x}_{k,j},y_{k,j}; \widetilde{{\bf w}}_k)}{M_k} -\frac{\ell'(\mathbf{x}_{k,M_k}',y_{k,M_k}'; \widetilde{{\bf w}}_k)}{M_k} - \frac{\lambda R'({\bf w}_k)}{K} \Bigr). 
\label{eqfirstprooflemmathree}
\end{align}

We notice that the primal updates corresponding with $\mathcal{D}$ and $\mathcal{D}'$ differ only for the $\ell$-update, where for the index $M_k$, the vector $\mathbf{x}_{k,M_k}$ and the scalar $y_{k,M_k}$ are different from $\mathbf{x}_{k,M_k}'$ and $y_{k,M_k}'$. Thus, for any neighboring data set $\mathcal{D}$ and $\mathcal{D}'$, the following holds:
\begin{align}
    &||{\bf w}_{k,\mathcal{D}}^{(n)} - {\bf w}_{k,\mathcal{D}'}^{(n)}|| = \Big|\Big| \frac{1}{2 \rho |\mathcal{N}_k| + \frac{1}{\eta^{(n)}}} \frac{1}{M_k}  \\
    &\Bigl( \ell'(\mathbf{x}_{k,M_k}, y_{k,M_k}, \widetilde{{\bf w}}_{k}^{(n-1)}) - \ell'(\mathbf{x}_{k,M_k}', y_{k,M_k}', \widetilde{{\bf w}}_{k}^{(n-1)}) \Bigr) \Big|\Big|. \notag
\end{align}
Since we assumed that $||\ell'(\cdot)||$ is bounded by $c_1$, the $l_2$-norm sensitivity is given by
\begin{align*}
    \max_{\mathcal{D},\mathcal{D}'} ||{\bf w}_{k,\mathcal{D}}^{(n)} - {\bf w}_{k,\mathcal{D}'}^{(n)}|| \leqslant \frac{2 c_1}{M_k (2 \rho |\mathcal{N}_k| + \frac{1}{\eta^{(n)}})}. \quad  \quad \qedhere
\end{align*}
\end{proof}

\section{Proof of Theorem I}
\begin{proof}
For any client $k$, at any step $t$, we add to the primal update a white Gaussian noise of variance $\sigma_{k}^{2}(n) \boldsymbol{I}_P$, that is equivalent to $\widetilde{{\bf w}}_{k}^{(n)} \sim \mathcal{N}({\bf w}_{k}^{(n)}, \sigma_{k}^{2}(n)\boldsymbol{I}_P)$. Hence, for two neighboring data sets $\mathcal{D}$ and $\mathcal{D}'$, we have $\widetilde{{\bf w}}_{k,\mathcal{D}}^{(n)} \sim \mathcal{N}({\bf w}_{k,\mathcal{D}}^{(n)}, \sigma_{k}^{2}(n)\boldsymbol{I}_P)$ and $\widetilde{{\bf w}}_{k,\mathcal{D}'}^{(n)} \sim \mathcal{N}({\bf w}_{k,\mathcal{D}'}^{(n)}, \sigma_{k}^{2}(n)\boldsymbol{I}_P)$.

Therefore, using \cite[Lemma 17]{Bun}, which states that $D_{\alpha}(N(\mu,\sigma^2\boldsymbol{I}_d) || N(\nu,\sigma^2\boldsymbol{I}_d)) = \frac{\alpha ||\mu - \nu ||_2^2}{2 \sigma^2}$, $\forall \alpha \in [1,\infty)$; we obtain, $\forall \alpha \in [1,\infty)$, the following R\'enyi divergence and simplification  using Lemma I:
\begin{align}
    D_{\alpha}(\widetilde{{\bf w}}_{k,\mathcal{D}}^{(n)} || \widetilde{{\bf w}}_{k,\mathcal{D}'}^{(n)}) = \frac{\alpha ||{\bf w}_{k,\mathcal{D}}^{(n)} - {\bf w}_{k,\mathcal{D}'}^{(n)}||_2^2}{2 \sigma_{k}^{2}(n)} \leqslant \frac{\alpha \Delta_k^2(n)}{2 \sigma_{k}^{2}(n)}.
\end{align}
We now consider the privacy loss of $\widetilde{{\bf w}}_{k}^{(n)}$ at output $\lambda$:
\begin{align}
   \boldsymbol{z}_k^{(n)}(\widetilde{{\bf w}}_{k,\mathcal{D}}^{(n)} || \widetilde{{\bf w}}_{k,\mathcal{D}'}^{(n)}) = \log \frac{P(\widetilde{{\bf w}}_{k,\mathcal{D}}^{(n)} = \lambda)}{P(\widetilde{{\bf w}}_{k,\mathcal{D}'}^{(n)} = \lambda)}.
\end{align}
As $D_{\alpha}(\cdot) \leqslant \epsilon + \rho \alpha \Longleftrightarrow E(e^{(\alpha-1)Z(\cdot)}) \leqslant e^{(\alpha - 1)(\epsilon + \rho \alpha)}$, we have:
\begin{align*}
    E(e^{(\alpha - 1)\boldsymbol{z}_k^{(n)}(\lambda)})
    \leqslant e^{(\alpha - 1) \frac{\alpha \Delta_k^2(n)}{2 \sigma_{k}^{2}(n)}}.
\end{align*}
Thus, the zCDP-NFL algorithm satisfies the dynamic $\varphi_{k}^{(n)}$-zCDP with $\varphi_{k}^{(n)} = \frac{\Delta_k^2(n)}{2 \sigma_{k}^{2}(n)}$.
\end{proof}

\section{Proof of Corollary}
\begin{proof}
Using \cite[Lemma 7]{Bun} and Theorem I, each client $k$ of the network has zCDP with $\varphi$ parameter $\sum_{0<n<T} \varphi_{k}^{(n)}$, $T$ being the final iteration index.

\noindent Since $\varphi_{k}^{(n+1)} = \varphi_{k}^{(n)}/\tau$, we have
\begin{align}
    \sum_{0<n<T} \varphi_{k}^{(n)} = \varphi_{k}^{(1)} \frac{1-\tau^{T}}{\tau^{T-1}-\tau^{T}}. \notag
\end{align}

Using \cite[Prop. 3]{Bun}, zCDP-NFL provides, $\forall \delta \in (0, 1)$, each client $k$ with $(\epsilon_k,\delta)$-DP, where $\epsilon_k = \varphi_{k}^{(1)} \frac{1-\tau^{T}}{\tau^{T-1}-\tau^{T}} + 2 \sqrt{\varphi_{k}^{(1)} \frac{1-\tau^{T}}{\tau^{T-1}-\tau^{T}} log \frac{1}{\delta}}$. Thus, the total privacy of the algorithm can be given in the DP metric with parameters $(\epsilon,\delta)$, $\forall \delta \in (0, 1)$, $\epsilon = \max_{k \in \mathcal{C}} \epsilon_k$.
\end{proof}

\section{Proof of Lemma II}
\begin{proof}

Using the convexity of $f(\cdot)$ we have:

\indent $f(\widetilde{\boldsymbol{\mathfrak{w}}}^{(n)}) - f(\boldsymbol{\mathfrak{w}}^{*}) \leqslant \langle \widetilde{\boldsymbol{\mathfrak{w}}}^{(n)} - \boldsymbol{\mathfrak{w}}^{*}, f'(\widetilde{\boldsymbol{\mathfrak{w}}}^{(n)})\rangle$.  And since $\hat{f}'(\widetilde{\boldsymbol{\mathfrak{w}}}^{(n+1)}) = 2 \boldsymbol{D}^{(n+1)} \otimes \boldsymbol{I}_P (\widetilde{\boldsymbol{\mathfrak{w}}}^{(n+1)} - \widetilde{\boldsymbol{\mathfrak{w}}}^{(n)}) + f'(\widetilde{\boldsymbol{\mathfrak{w}}}^{(n)})$, we have $f'(\widetilde{\boldsymbol{\mathfrak{w}}}^{(n)}) = \hat{f}'(\widetilde{\boldsymbol{\mathfrak{w}}}^{(n+1)}) -  2 \boldsymbol{D}^{(n+1)} \otimes \boldsymbol{I}_P (\widetilde{\boldsymbol{\mathfrak{w}}}^{(n+1)} - \widetilde{\boldsymbol{\mathfrak{w}}}^{(n)})$.

Combining both equations we obtain:
\begin{align}
\label{ComposedConvexity}
    &f(\widetilde{\boldsymbol{\mathfrak{w}}}^{(n)}) - f(\boldsymbol{\mathfrak{w}}^{*}) \leqslant \langle \widetilde{\boldsymbol{\mathfrak{w}}}^{(n)} - \boldsymbol{\mathfrak{w}}^{*}, \hat{f}'(\widetilde{\boldsymbol{\mathfrak{w}}}^{(n+1)}) - 2 \boldsymbol{D}^{(n+1)} \otimes \boldsymbol{I}_P (\widetilde{\boldsymbol{\mathfrak{w}}}^{(n+1)} - \widetilde{\boldsymbol{\mathfrak{w}}}^{(n)})\rangle.
\end{align}

Employing \eqref{CVproofIntroLast} in \eqref{ComposedConvexity}, followed by some algebraic manipulations, yields


\begin{align}
\label{LastBeforeQ}
    &\frac{f(\widetilde{\boldsymbol{\mathfrak{w}}}^{(n)}) - f(\boldsymbol{\mathfrak{w}}^{*})}{\rho} + \langle\widetilde{\boldsymbol{\mathfrak{w}}}^{(n)}, 2 \boldsymbol{Q} \boldsymbol{r}\rangle \notag \\
    &\leqslant \langle\widetilde{\boldsymbol{\mathfrak{w}}}^{(n)}, 2 \boldsymbol{Q} \boldsymbol{r}\rangle+  \langle \widetilde{\boldsymbol{\mathfrak{w}}}^{(n)} - \boldsymbol{\mathfrak{w}}^{*}, 2 \boldsymbol{M} \boldsymbol{\xi}^{(n+1)} - 2 \boldsymbol{Q} \boldsymbol{r}^{(n+1)} \notag\\
    & \quad - \boldsymbol{L}_+(\widetilde{\boldsymbol{\mathfrak{w}}}^{(n+1)} - \widetilde{\boldsymbol{\mathfrak{w}}}^{(n)})    - \frac{2}{\rho} \boldsymbol{D}^{(n+1)} \otimes \boldsymbol{I}_P (\widetilde{\boldsymbol{\mathfrak{w}}}^{(n+1)} - \widetilde{\boldsymbol{\mathfrak{w}}}^{(n)})\rangle, \notag\\
    &\leqslant   \langle \widetilde{\boldsymbol{\mathfrak{w}}}^{(n)} - \boldsymbol{\mathfrak{w}}^{*},   2 \boldsymbol{Q} (\boldsymbol{r} - \boldsymbol{r}^{(n+1)} )  + \boldsymbol{L}_+(\widetilde{\boldsymbol{\mathfrak{w}}}^{(n)} - {\boldsymbol{\mathfrak{w}}}^{(n+1)}) \notag\\
    & \quad + \boldsymbol{L}_-(\widetilde{\boldsymbol{\mathfrak{w}}}^{(n+1)} - {\boldsymbol{\mathfrak{w}}}^{(n+1)})  - \frac{2}{\rho} \boldsymbol{D}^{(n+1)} \otimes \boldsymbol{I}_P (\widetilde{\boldsymbol{\mathfrak{w}}}^{(n+1)} - \widetilde{\boldsymbol{\mathfrak{w}}}^{(n)})\rangle. 
\end{align}

It follows that
\begin{align}
    || \boldsymbol{q}^{(n)} - \boldsymbol{q}^*||_{\boldsymbol{G}}^2 &= \left\langle \begin{pmatrix} \boldsymbol{r}^{(n)} - \boldsymbol{r}\\ \widetilde{\boldsymbol{\mathfrak{w}}}^{(n)} - \boldsymbol{\mathfrak{w}}^{*} \end{pmatrix}, \begin{pmatrix} \rho( \boldsymbol{r}^{(n)} - \boldsymbol{r})\\ \frac{\rho \boldsymbol{L}_+}{2} (\widetilde{\boldsymbol{\mathfrak{w}}}^{(n)} - \boldsymbol{\mathfrak{w}}^{*}) \end{pmatrix} \right\rangle \notag\\
    &= \rho || \boldsymbol{r}^{(n)} - \boldsymbol{r}||_2^2 + ||\widetilde{\boldsymbol{\mathfrak{w}}}^{(n)} - \boldsymbol{\mathfrak{w}}^{*}||_{\frac{\rho \boldsymbol{L}_+}{2}}^2. \notag
\end{align}

In particular, we obtain the equality:
\begin{align}
\label{Qeq}
   &\frac{1}{\rho}(||\boldsymbol{q}^{(n-1)} - \boldsymbol{q}^{*}||^2_{\boldsymbol{G}} - ||\boldsymbol{q}^{(n)} - \boldsymbol{q}^{*}||^2_{\boldsymbol{G}} - ||\boldsymbol{q}^{(n)} - \boldsymbol{q}^{(n-1)}||^2_{\boldsymbol{G}}) \notag\\
   &= \langle \widetilde{\boldsymbol{\mathfrak{w}}}^{(n)} - \boldsymbol{\mathfrak{w}}^{*}, 2 \boldsymbol{Q} (\boldsymbol{r} - \boldsymbol{r}^{(n)})\rangle + \langle \widetilde{\boldsymbol{\mathfrak{w}}}^{(n)} - \boldsymbol{\mathfrak{w}}^{*}, \boldsymbol{L}_+(\widetilde{\boldsymbol{\mathfrak{w}}}^{(n-1)} - \widetilde{\boldsymbol{\mathfrak{w}}}^{(n)}) \rangle.
\end{align}

Which we use to reformulate the second term of \eqref{LastBeforeQ}, after which it is combined with \eqref{Qeq}. Using the fact that $ \boldsymbol{Q} \widetilde{\boldsymbol{\mathfrak{w}}}^{(n+1)} = \boldsymbol{r}^{(n+1)} - \boldsymbol{r}^{(n)}$, we obtain:

\begin{align}&\frac{f(\widetilde{\boldsymbol{\mathfrak{w}}}^{(n)}) - f(\boldsymbol{\mathfrak{w}}^{*})}{\rho} + \langle\widetilde{\boldsymbol{\mathfrak{w}}}^{(n)}, 2 \boldsymbol{Q} \boldsymbol{r}\rangle 
    \leqslant \frac{1}{\rho}(||\boldsymbol{q}^{(n-1)} - \boldsymbol{q}^{*}||^2_{\boldsymbol{G}} - ||\boldsymbol{q}^{(n)} - \boldsymbol{q}^{*}||^2_{\boldsymbol{G}}) \notag\\
    &\quad - \frac{\Phi_{\min}(\boldsymbol{L}_-)}{2} || \widetilde{\boldsymbol{\mathfrak{w}}}^{(n)} - \boldsymbol{\mathfrak{w}}^* ||^2_2 - ||\widetilde{\boldsymbol{\mathfrak{w}}}^{(n)} - \widetilde{\boldsymbol{\mathfrak{w}}}^{(n-1)}||_{\frac{\boldsymbol{L}_+}{2}}^2   - 2 \langle \boldsymbol{Q} \widetilde{\boldsymbol{\mathfrak{w}}}^{(n)}, \boldsymbol{Q} \widetilde{\boldsymbol{\mathfrak{w}}}^{(n+1)} \rangle \notag\\
    &\quad + \langle \widetilde{\boldsymbol{\mathfrak{w}}}^{(n)} - \boldsymbol{\mathfrak{w}}^{*}, \boldsymbol{L}_+(2 \widetilde{\boldsymbol{\mathfrak{w}}}^{(n)} - \widetilde{\boldsymbol{\mathfrak{w}}}^{(n-1)} - \widetilde{\boldsymbol{\mathfrak{w}}}^{(n+1)}) \rangle \notag\\
    &\quad + ||\widetilde{\boldsymbol{\mathfrak{w}}}^{(n)} - \boldsymbol{\mathfrak{w}}^{*}||_2 ||\boldsymbol{L}_+ \boldsymbol{\xi}^{(n+1)}||_2   + ||\widetilde{\boldsymbol{\mathfrak{w}}}^{(n)} - \boldsymbol{\mathfrak{w}}^{*}||_2 ||\boldsymbol{L}_- \boldsymbol{\xi}^{(n+1)}||_2 \notag\\
    &\quad + \langle \widetilde{\boldsymbol{\mathfrak{w}}}^{(n)} - \boldsymbol{\mathfrak{w}}^{*}, \frac{2}{\rho} \boldsymbol{D}^{(n+1)} \otimes \boldsymbol{I}_P (\widetilde{\boldsymbol{\mathfrak{w}}}^{(n)} - \widetilde{\boldsymbol{\mathfrak{w}}}^{(n+1)})\rangle \notag
\end{align}

Using the inequality $||\boldsymbol{a}|| ||\boldsymbol{b}|| \leqslant m ||\boldsymbol{a}|| + \frac{1}{m} ||\boldsymbol{b}||$ for $m > 0$ with $m=\Phi_{\min}(\boldsymbol{L}_-)$, it can be established that
\begin{align}
    &\frac{f(\widetilde{\boldsymbol{\mathfrak{w}}}^{(n)}) - f(\boldsymbol{\mathfrak{w}}^{*})}{\rho} + \langle\widetilde{\boldsymbol{\mathfrak{w}}}^{(n)}, 2 \boldsymbol{Q} \boldsymbol{r}\rangle  
     \leqslant \frac{1}{\rho}(||\boldsymbol{q}^{(n-1)} - \boldsymbol{q}^{*}||^2_{\boldsymbol{G}} - ||\boldsymbol{q}^{(n)} - \boldsymbol{q}^{*}||^2_{\boldsymbol{G}}) \notag\\
    &\quad - 2 \langle \boldsymbol{Q} \widetilde{\boldsymbol{\mathfrak{w}}}^{(n)}, \boldsymbol{Q} \widetilde{\boldsymbol{\mathfrak{w}}}^{(n+1)} \rangle - ||\widetilde{\boldsymbol{\mathfrak{w}}}^{(n)} - \widetilde{\boldsymbol{\mathfrak{w}}}^{(n-1)}||_{\frac{\boldsymbol{L}_+}{2}}^2 \notag\\
    &\quad + \langle \widetilde{\boldsymbol{\mathfrak{w}}}^{(n)} - \boldsymbol{\mathfrak{w}}^{*}, \boldsymbol{L}_+(2 \widetilde{\boldsymbol{\mathfrak{w}}}^{(n)} - \widetilde{\boldsymbol{\mathfrak{w}}}^{(n-1)} - \widetilde{\boldsymbol{\mathfrak{w}}}^{(n+1)}) \rangle \notag\\
    &\quad + \frac{4(\Phi_{\max}(\boldsymbol{L}_-)^2 + \Phi_{\max}(\boldsymbol{L}_+)^2)}{\Phi_{\min}(\boldsymbol{L}_-)} ||\boldsymbol{\xi}^{(n+1)}||^2_2 \notag\\
    &\quad + \langle \widetilde{\boldsymbol{\mathfrak{w}}}^{(n)} - \boldsymbol{\mathfrak{w}}^{*}, \frac{2}{\rho} \boldsymbol{D}^{(n+1)} \otimes \boldsymbol{I}_P (\widetilde{\boldsymbol{\mathfrak{w}}}^{(n)} - \widetilde{\boldsymbol{\mathfrak{w}}}^{(n+1)})\rangle. \notag \qedhere
\end{align}


\end{proof}

\section{Proof of Theorem II}

\begin{proof}

We first take the sum of the result of Lemma II from $n=1$ to $n=T$ to obtain a bound given by
\begin{align}
    \label{StartThII}
    &\frac{1}{\rho} (\sum_{n=1}^{T} f(\widetilde{\boldsymbol{\mathfrak{w}}}^{(n)}) - f(\boldsymbol{\mathfrak{w}}^*)) + \langle2 \boldsymbol{r}, \sum_{n=1}^{T} \boldsymbol{Q} \widetilde{\boldsymbol{\mathfrak{w}}}^{(n)}\rangle \\
    & \leqslant \sum_{n=1}^{T} \Bigl( \frac{4(\Phi_{\max}(\boldsymbol{L}_-)^2 + \Phi_{\max}(\boldsymbol{L}_+)^2)}{\Phi_{\min}(\boldsymbol{L}_-)} ||\boldsymbol{\xi}^{(n+1)}||^2_2 - 2 \langle \boldsymbol{Q} \widetilde{\boldsymbol{\mathfrak{w}}}^{(n)}, \boldsymbol{Q} \widetilde{\boldsymbol{\mathfrak{w}}}^{(n+1)} \rangle\notag\\
    &\qquad   - ||\widetilde{\boldsymbol{\mathfrak{w}}}^{(n)} - \widetilde{\boldsymbol{\mathfrak{w}}}^{(n-1)}||_{\frac{\boldsymbol{L}_+}{2}}^2   + \langle \widetilde{\boldsymbol{\mathfrak{w}}}^{(n)} - \boldsymbol{\mathfrak{w}}^{*}, \frac{2}{\rho} \boldsymbol{D}^{(n+1)} \otimes \boldsymbol{I}_P (\widetilde{\boldsymbol{\mathfrak{w}}}^{(n)} - \widetilde{\boldsymbol{\mathfrak{w}}}^{(n+1)})\rangle \notag\\
    &\qquad + \langle \widetilde{\boldsymbol{\mathfrak{w}}}^{(n)} - \boldsymbol{\mathfrak{w}}^{*}, \boldsymbol{L}_+(2 \widetilde{\boldsymbol{\mathfrak{w}}}^{(n)} - \widetilde{\boldsymbol{\mathfrak{w}}}^{(n-1)} - \widetilde{\boldsymbol{\mathfrak{w}}}^{(n+1)}) \rangle \Bigr)  + \frac{1}{\rho} || \boldsymbol{q}^{(0)} - \boldsymbol{q}^*||_{\boldsymbol{G}}^2. \notag
\end{align}

\noindent After setting $r=0$ and performing some algebraic manipulations, we obtain:
\begin{align}
    &\frac{1}{\rho} (\sum_{n=1}^{T} f(\widetilde{\boldsymbol{\mathfrak{w}}}^{(n+1)}) - f(\boldsymbol{\mathfrak{w}}^*))  \leqslant \sum_{n=1}^{T} \Bigl( \frac{4(\Phi_{\max}(\boldsymbol{L}_-)^2 + \Phi_{\max}(\boldsymbol{L}_+)^2)}{\Phi_{\min}(\boldsymbol{L}_-)} ||\boldsymbol{\xi}^{(n+1)}||^2_2 \notag\\
    &\qquad - 2 \langle \boldsymbol{Q} \widetilde{\boldsymbol{\mathfrak{w}}}^{(n)}, \boldsymbol{Q} \widetilde{\boldsymbol{\mathfrak{w}}}^{(n+1)} \rangle  - ||\widetilde{\boldsymbol{\mathfrak{w}}}^{(n)} - \widetilde{\boldsymbol{\mathfrak{w}}}^{(n-1)}||_{\frac{\boldsymbol{L}_+}{2}}^2 \notag\\
    &\qquad + \langle \widetilde{\boldsymbol{\mathfrak{w}}}^{(n)} - \boldsymbol{\mathfrak{w}}^{*}, \frac{2}{\rho} \boldsymbol{D}^{(n+1)} \otimes \boldsymbol{I}_P (\widetilde{\boldsymbol{\mathfrak{w}}}^{(n)} - \widetilde{\boldsymbol{\mathfrak{w}}}^{(n+1)})\rangle \notag\\
    &\qquad - \langle\boldsymbol{\mathfrak{w}}^{*}, \boldsymbol{L}_+(2 \widetilde{\boldsymbol{\mathfrak{w}}}^{(n)} - \widetilde{\boldsymbol{\mathfrak{w}}}^{(n-1)} - \widetilde{\boldsymbol{\mathfrak{w}}}^{(n+1)}) \rangle    + ||\widetilde{\boldsymbol{\mathfrak{w}}}^{(n+1)} - \widetilde{\boldsymbol{\mathfrak{w}}}^{(n)}||_{\boldsymbol{L}_+}^2 \Bigr) \notag\\
    &\quad + \langle \widetilde{\boldsymbol{\mathfrak{w}}}^{(1)}, \boldsymbol{L}_+(\widetilde{\boldsymbol{\mathfrak{w}}}^{(1)} - \widetilde{\boldsymbol{\mathfrak{w}}}^{(0)})\rangle   + || \boldsymbol{Q} \widetilde{\boldsymbol{\mathfrak{w}}}^{(0)}||_2^2 + ||\widetilde{\boldsymbol{\mathfrak{w}}}^{(0)} - \boldsymbol{\mathfrak{w}}^*||_{\frac{ \boldsymbol{L}_-}{2}}^2.     \label{BeforeExpThII}
\end{align}
Using Jensen's inequality on the expectation of \eqref{BeforeExpThII}, we obtain:
\begin{align}
    &\mathbb{E}[f(\hat{\boldsymbol{\mathfrak{w}}}^{(T)}) - f(\boldsymbol{\mathfrak{w}}^*)] \leqslant \frac{\rho}{T} \sum_{n=1}^{T} \Bigl( - 2 \langle \boldsymbol{Q} \widetilde{\boldsymbol{\mathfrak{w}}}^{(n)}, \boldsymbol{Q} \widetilde{\boldsymbol{\mathfrak{w}}}^{(n+1)} \rangle - ||\widetilde{\boldsymbol{\mathfrak{w}}}^{(n)} - \widetilde{\boldsymbol{\mathfrak{w}}}^{(n-1)}||_{\frac{\boldsymbol{L}_+}{2}}^2 \notag\\
    &\qquad + \langle \widetilde{\boldsymbol{\mathfrak{w}}}^{(n)} - \boldsymbol{\mathfrak{w}}^{*}, \frac{2}{\rho} \boldsymbol{D}^{(n+1)} \otimes \boldsymbol{I}_P (\widetilde{\boldsymbol{\mathfrak{w}}}^{(n)} - \widetilde{\boldsymbol{\mathfrak{w}}}^{(n+1)})\rangle \notag\\
    &\qquad - \langle\boldsymbol{\mathfrak{w}}^{*}, \boldsymbol{L}_+(2 \widetilde{\boldsymbol{\mathfrak{w}}}^{(n)} - \widetilde{\boldsymbol{\mathfrak{w}}}^{(n-1)} - \widetilde{\boldsymbol{\mathfrak{w}}}^{(n+1)}) \rangle \notag\\
    &\qquad + ||\widetilde{\boldsymbol{\mathfrak{w}}}^{(n+1)} - \widetilde{\boldsymbol{\mathfrak{w}}}^{(n)}||_{\boldsymbol{L}_+}^2 \Bigr)  + \frac{\langle \widetilde{\boldsymbol{\mathfrak{w}}}^{(1)}, \boldsymbol{L}_+(\widetilde{\boldsymbol{\mathfrak{w}}}^{(1)} - \widetilde{\boldsymbol{\mathfrak{w}}}^{(0)})\rangle}{T} + \frac{\rho || \boldsymbol{Q} \widetilde{\boldsymbol{\mathfrak{w}}}^{(0)}||_2^2}{T}  \notag\\
    &\quad + \frac{1}{T} \frac{\rho P 4(\Phi_{\max}(\boldsymbol{L}_-)^2 + \Phi_{\max}(\boldsymbol{L}_+)^2) \sum_{k=1}^K \sigma^{2 (0)}_k} {\Phi_{\min}(\boldsymbol{L}_-) (1 - \tau)}  + \frac{ \rho ||\widetilde{\boldsymbol{\mathfrak{w}}}^{(0)} - \boldsymbol{\mathfrak{w}}^*||_{\frac{ \boldsymbol{L}_-}{2}}^2}{T}. \notag \qedhere
\end{align}
\end{proof}

 \bibliographystyle{elsarticle-num} 
 \bibliography{References}





\end{document}